\documentclass[final,hidelinks,onefignum,onetabnum]{siamart220329}
\usepackage{amsfonts, amsmath, amssymb} 
\usepackage{mathtools, mathrsfs}
\usepackage{stmaryrd}
\usepackage[algo2e,linesnumbered,commentsnumbered,ruled,vlined]{algorithm2e}
\usepackage[english]{babel}
\usepackage{pgf}
\usepackage{graphicx}
\usepackage{caption}
\usepackage{subcaption}
\usepackage{pgfplots, pgfplotstable}
\pgfplotsset{compat=1.5}
\usetikzlibrary{spy}
\ifpdf
  \DeclareGraphicsExtensions{.eps,.pdf,.png,.jpg}
\else
  \DeclareGraphicsExtensions{.eps}
\fi

\newcommand{\rev}[1]{{#1}}

\newcommand\ONECOLWIDTH{0.95\textwidth}
\newcommand\TWOCOLWIDTH{0.40\textwidth}

\newcommand{\mute}[1]{}

\newsiamremark{remark}{Remark}
\newsiamremark{hypothesis}{Hypothesis}
\crefname{hypothesis}{Hypothesis}{Hypotheses}
\newsiamthm{claim}{Claim}

\newcommand{\RR}{\mathbb{R}}

\DeclarePairedDelimiter\abs{\lvert}{\rvert}

\DeclarePairedDelimiter\jump{\llbracket}{\rrbracket}

\DeclareMathOperator{\dist}{dist}

\newcommand{\om}{\Omega}
\newcommand{\oh}{\Omega_h}
\newcommand{\veps}{\varepsilon}
\newcommand{\Ho}{\mathring{H}}
\newcommand{\ut}{\tilde{u}}
\newcommand{\vt}{\tilde{v}}
\newcommand{\rhot}{\tilde{\rho}}
\newcommand{\Eth}{E_{\tau, h}}
\newcommand{\tEth}{\tilde{E}_{\tau, h}}

\DeclareMathOperator{\diver}{div}
\DeclareMathOperator{\grad}{grad}
\DeclareMathOperator{\curl}{curl}
\DeclareMathOperator{\rot}{rot}
\DeclareMathOperator{\dive}{div}

\headers{Adaptive resolution of fine scales in modes of MOF}
{J.~Gopalakrishnan, J.~Grosek, G.~Pinochet-Soto, and P.~Vandenberge}

\title{Adaptive resolution of fine scales in modes of microstructured optical fibers
    \thanks{Submitted to the editors 04/04/2024.
    \funding{This work was supported in part by AFOSR Grant FA9550-23-1-0103 and NSF Grant 2245077. The work also benefited from activities organized under the auspices of NSF RTG Grant 2136228. }}
}

\author{
    Jay Gopalakrishnan
    \thanks{
    The Fariborz Maseeh Department of Mathematics + Statistics, Portland State University,
    PO Box 751, Portland, OR, USA
    (\email{gjay@pdx.edu}, \email{gpin2@pdx.edu}, \email{piet2@pdx.edu}).}
    \and
    Jacob Grosek
    \thanks{
    Directed Energy Directorate, Air Force Research Laboratory,
    3550 Aberdeen Ave SE, Kirtland Air Force Base, Albuquerque, NM (\email{AFRL.RDL.OrgBox@us.af.mil}).}
    \and
    Gabriel Pinochet-Soto\footnotemark[2]
    \and
    Pieter Vandenberge\footnotemark[2]
}

\usepackage{amsopn}

\ifpdf
\hypersetup{
  pdftitle={Adaptive resolution of fine scales in modes of microstructured optical fibers},
  pdfauthor={J. Gopalakrishnan, J. Grosek, G. Pinochet-Soto, and P. Vandenberge}
}
\fi

\begin{document}

\maketitle

\begin{abstract}
	An adaptive algorithm for computing eigenmodes and propagation constants of optical fibers is proposed. 
	The algorithm is built using a dual-weighted residual error estimator. 
	The residuals are based on the eigensystem for leaky hybrid modes obtained from Maxwell equations truncated to a finite domain after a transformation by a perfectly matched layer. 
	The adaptive algorithm is then applied to compute practically interesting modes for multiple fiber microstructures. 
	Emerging microstructured optical fibers are characterized by complex geometrical features in their transverse cross-section. 
	Their leaky modes, useful for confining and propagating light in their cores, often exhibit fine scale features. 
	The adaptive algorithm automatically captures these features without any expert input. 
	The results also show that confinement losses of these modes are captured accurately on the adaptively found meshes. 
\end{abstract}



\section{Introduction}\label{sec:intro}

Microstructured optical fibers with complex features are increasingly being fabricated and used~\cite{CordeNgEbend20}. 
An important consideration in their design is their confinement loss or radiation loss. 
This loss can be extracted from a complex eigenvalue of an eigenproblem arising from the Maxwell system governing the propagation of light within optical fibers. 
The corresponding eigenfunctions are practically interesting since they form the leaky modes of fibers capable of propagating substantial portions of the input energy through the fiber. 
The electric fields of these modes generally have both transverse and longitudinal vector components, so they are also referred to as ``hybrid modes'' \cite{SnydeLove83}. 
The purpose of this paper is to give an adaptive algorithm for capturing such modes with sufficient accuracy. 
The algorithm is based on the mathematical technique of dual-weighted residual (DWR) error estimators~\cite{2001BecRan}. 
We detail how to apply this approach to a finite element discretization of the eigenproblem for hybrid modes combined with a domain truncation after a transformation by the perfectly matched layer (PML) \cite{1994Ber, ColliMonk98}. 

Not all microstructured fiber designs use the same physical mechanism for guiding light. 
Some fibers use confinement by anti-resonant reflection effects \cite{2014Pol, DuguaKokubKoch86, 2013KolKosPryBirPlotDia} while others use confinement by resonances of guided modes, total internal reflection, and photonic confinement within defects of periodic arrangements~\cite{GuKongHawki15,Knigh03, OhterHirosYamad14}. 
The diversity of these physical mechanisms necessitate expert guidance in each fiber's numerical simulation, thus posing challenges in creating general purpose simulation tools for these fibers. 
This paper contributes to addressing this challenge. 
Specifically, several practically useful leaky modes of varied microstructured fibers have been found to exhibit fine-scale ripples in varied locations. 
As shown in~\cite{2023VanGopGro}, accuracy of computed confinement losses depend on resolving such fine-scale features. 
To do so without expert insight and to facilitate design optimization, an algorithm that can automatically detect such fine-scale features accompanied by appropriate mesh refinement is useful. 
This paper provides such an algorithm, basing it on sound mathematical principles independent of fiber design. 
A prior work proposing a goal-oriented adaptivity strategy for computing losses accurately is~\cite{2007ZscBurPom}, but as best as we can see, the goal functional proposed there is not a continuous functional and we do not know how to make their technique mathematically rigorous. 
Instead, we develop an alternate approach based on the ideas in~\cite{2001HeuRan} which directly target the eigenvalue error. 

The eigenvalues of the Maxwell system for hybrid modes are squares of the physical propagation constant~$\beta$ of the modes. 
Confinement loss of a fiber is a scalar multiple of the imaginary part of $\beta$. 
Typical values of $\beta$ have real parts of the order of $10^6$ and imaginary parts, $\text{Im}(\beta)$, that can range in orders of magnitude $10^{-6}$ to $10^{3}$, depending on the fiber and operating wavelength.
To get around this large variation in order of magnitudes, instead of a mere length scaling, we perform a nondimensionalization of the Maxwell system (see Section~\ref{sec:pmlmaxwell}) akin to the process of going from Helmholtz to a Schr\"odinger problem. 
Then, instead of $\beta^2$, we compute a nondimensional $Z^2$, of unit order of magnitude, in the complex plane. 
Applying this approach to varied fibers, we find that an adaptive algorithm tailored to reduce the error in $Z^2$ also appears to get the confinement loss accurately and capture the above-mentioned fine-scale ripples. Whether a mathematically rigorous algorithm tailored specifically to reduce the error in loss can be found, and if so, whether it will do better, is an issue for further research. 

The basis for DWR-based adaptive algorithms is developed in terms of
optimal control in the review~\cite{2001BecRan}. Roughly speaking,
their idea is that the error in a quantity of interest can be
estimated by an element-wise residual (which treats local error
contributions), weighted by the solution of a dual problem (which
takes global error transport into account).  DWR algorithms specific
for eigenvalues have been outlined for elliptic problems
in~\cite{2001HeuRan,RannaWesteWolln10}. \rev{We choose the DWR error
  estimation technique because our eigenproblem is nonselfadjoint and
  because our eigensolver readily gives the dual solution (left
  eigenvector) in addition to the solution (the eigenvalue and the
  right eigenvector). Not incurring extra costs for solving the dual
  problem in this case makes the DWR technique an efficient choice.
To apply the DWR methodology to our Maxwell eigenvalue application,}   
we bring in an additional technique, namely the use of a regular
decomposition. The resulting 
theory is summarized in a single self-contained theorem and
proof in Section~\ref{sec:dwr-maxwell} giving the error estimator
expression, after which we describe our adaptive algorithm.  Later
sections clearly reveal the efficacy of the algorithm for computing
leaky modes of complex fibers.  In particular, Section~\ref{sec:bragg}
is devoted to verifying the correctness of the methodology using a
Bragg fiber for which mode solutions in exact closed form are available.

\begin{figure}[htbp]
	\centering
	\begin{subfigure}[b]{\TWOCOLWIDTH}
		\centering
		\begin{tikzpicture}[scale=0.035]
			\fill [blue!30,even odd rule] (0,0)
			circle[radius=70.8]
			circle[radius=38.0];
			\draw[red, <->] (0,0) -- (38.0,0)
			node [scale=0.75,black, midway, below]
				{$r_{\text{core}}$} ;
			\draw[red, <->] (-70.8,0) -- (-38.0,0)
			node [scale=0.75,black, midway, above]
				{$t_{\text{ring}}$} ;
			\draw[red, dashed] (0,0) circle [radius=82.8];
			\draw[red, <->] (0,0) -- (43.0,72.0)
			node [scale=0.75,black, midway, left]
				{$r_0$} ;
		\end{tikzpicture}
		\caption{Geometry of the Bragg fiber.}
		\label{fig:geo_bragg}
	\end{subfigure}
	\hfill
	\begin{subfigure}[b]{\TWOCOLWIDTH}
		\centering
		\begin{tikzpicture}[scale=0.035]
			\fill [blue!30,even odd rule] (0,0)
			circle[radius=70.8]
			circle[radius=38.0];
			\fill [blue!30,even odd rule] (0.0, 26.505) 
			circle [radius=12]
			circle [radius=10];
			\fill [blue!30,even odd rule] (-22.954003327306545, 13.2525) 
			circle [radius=12]
			circle [radius=10];
			\fill [blue!30,even odd rule] (-22.954003327306545, -13.2525) 
			circle [radius=12]
			circle [radius=10];
			\fill [blue!30,even odd rule] (0.0, -26.505) 
			circle [radius=12]
			circle [radius=10];
			\fill [blue!30,even odd rule] (22.954003327306545, -13.2525) 
			circle [radius=12]
			circle [radius=10];
			\fill [blue!30,even odd rule] (22.954003327306545, 13.2525) 
			circle [radius=12]
			circle [radius=10];
			\draw[red, dashed] (0,0) circle [radius=14.0];
			\draw[red, <->] (0,0.0) -- (14.0,0.0)
			node [scale=0.75,black, midway, below] 
				{$r_{\text{core}}$};
			\draw[red, |->] (5.0, -26.505) -- (10.0, -26.505);
			\draw[red, <-|] (12.0, -26.505) -- (17.0, -26.505);
			\node[below, scale=0.75] at (11.0, -26.505) {$t_{\text{cap}}$};
			\draw[red, <->] (-70.8,0.0) -- (-38.0,0.0)
			node [scale=0.75,black, midway, above] 
				{$t_{\text{clad}}$};
			\draw[red, dashed] (0,0) circle [radius=82.8];
			\draw[red, <->] (0,0) -- (43.0,72.0)
			node [scale=0.75,black, midway, left]
				{$r_0$} ;
			\draw[red, |->] (-20, -16) -- (-13.0, -20.0);
			\draw[red, <-|] (-10.75, -21.25) -- (-3.0, -26.0);
			\node[below, scale=0.75] at (-15.0, -16.5) {$d$};
		\end{tikzpicture}
		\caption{Geometry of the ARF fiber.}
		\label{fig:geo_arf}
	\end{subfigure}
	\vspace{1em}
	\begin{subfigure}[b]{\TWOCOLWIDTH}
		\centering
		\begin{tikzpicture}[scale = 0.035]
			\fill [blue!30,even odd rule] (0,0)
			circle[radius=70.8]
			circle[radius=38.0];
			\fill [blue!30,even odd rule] (0.0, 26.505) 
			circle [radius=12]
			circle [radius=10];
			\fill [blue!30,even odd rule] (0.0, 32.505) 
			circle [radius=6]
			circle [radius=4];
			\fill [blue!30,even odd rule] (0.0, -26.505) 
			circle [radius=12]
			circle [radius=10];
			\fill [blue!30,even odd rule] (0.0, -32.505) 
			circle [radius=6]
			circle [radius=4];
			\fill [blue!30,even odd rule] (-22.954003327306545, 13.2525) 
			circle [radius=12]
			circle [radius=10];
			\fill [blue!30,even odd rule] (-27.854003327306545, 16.8525)
			circle [radius=6]
			circle [radius=4];
			\fill [blue!30,even odd rule] (-22.954003327306545, -13.2525) 
			circle [radius=12]
			circle [radius=10];
			\fill [blue!30,even odd rule] (-27.854003327306545, -16.8525)
			circle [radius=6]
			circle [radius=4];
			\fill [blue!30,even odd rule] (22.954003327306545, -13.2525) 
			circle [radius=12]
			circle [radius=10];
			\fill [blue!30,even odd rule] (27.854003327306545, -16.8525)
			circle [radius=6]
			circle [radius=4];
			\fill [blue!30,even odd rule] (22.954003327306545, 13.2525) 
			circle [radius=12]
			circle [radius=10];
			\fill [blue!30,even odd rule] (27.854003327306545, 16.8525)
			circle [radius=6]
			circle [radius=4];
			\draw[red, |->] (5.0, -26.505) -- (10.0, -26.505);
			\draw[red, <-|] (12.0, -26.505) -- (17.0, -26.505);
			\node[below, scale=0.75] at (11.0, -26.505) {$t_{\text{cap, o}}$};
			\draw[red, dashed] (0,0) circle [radius=14.0];
			\draw[red, <->] (0,0.0) -- (14.0,0.0)
			node [scale=0.75,black, midway, below] 
				{$r_{\text{core}}$};
			\draw[red, <->] (-70.8,0.0) -- (-38.0,0.0)
			node [scale=0.75,black, midway, above] 
				{$t_{\text{clad}}$};
            	\draw[red, dashed] (0,0) circle [radius=82.8];
            	\draw[red, <->] (0,0) -- (42.0,72.0)
            	node [scale=0.75,black, midway, left]{$r_0$} ;
			\draw[red, |->] (-20, -16) -- (-13.0, -20.0);
			\draw[red, <-|] (-10.75, -21.25) -- (-3.0, -26.0);
			\node[below, scale=0.75] at (-15.0, -16.5) {$d$};
		\end{tikzpicture}
		\caption{Geometry of the NANF fiber.}
		\label{fig:geo_nanf}
	\end{subfigure}
	\hfill
	\begin{subfigure}[b]{\TWOCOLWIDTH}    
		\centering
		\begin{tikzpicture}[scale=0.035]
			\fill [blue!30] (0,0)
			circle[radius=70.8];
			\fill [blue!50] (26.505, 0.0)
			circle [radius=12];
			\fill [blue!50] (-26.505, 0.0)
			circle [radius=12];
			\fill [blue!50] (13.2525, 22.954003327306545)
			circle [radius=12];
			\fill [blue!50] (-13.2525, 22.954003327306545)
			circle [radius=12];
			\fill [blue!50] (-13.2525, -22.954003327306545)
			circle [radius=12];
			\fill [blue!50] (13.2525, -22.954003327306545)
			circle [radius=12];
			\draw[red, dashed] (0,0) circle [radius=14.0];
			\draw[red, <->] (0,0.0) -- (14.0,0.0)
				node [scale=0.75,black, midway, below] 
				{$r_{\text{core}}$};
			\draw[red, <->] (-70.8,0.0) -- (0.0,0.0)
				node [scale=0.75,black, midway, above] 
				{$r_{\text{outer}}$};
			\draw[red, <->] (1.25, 23.005) -- (13.25, 23.005)
				node [scale=0.75,black, midway, above] 
				{$r_{\text{tube}}$};
		\end{tikzpicture}
		\caption{Geometry of the PBG fiber.}
		\label{fig:geo_pbg}
	\end{subfigure}
	\vspace{1em}
	\caption{Transverse geometries of the microstructured fibers (not to scale) studied in this work. White regions indicate air and shades of blue indicate dielectric materials with higher refractive indices.}
	\label{fig:geometries}
\end{figure}
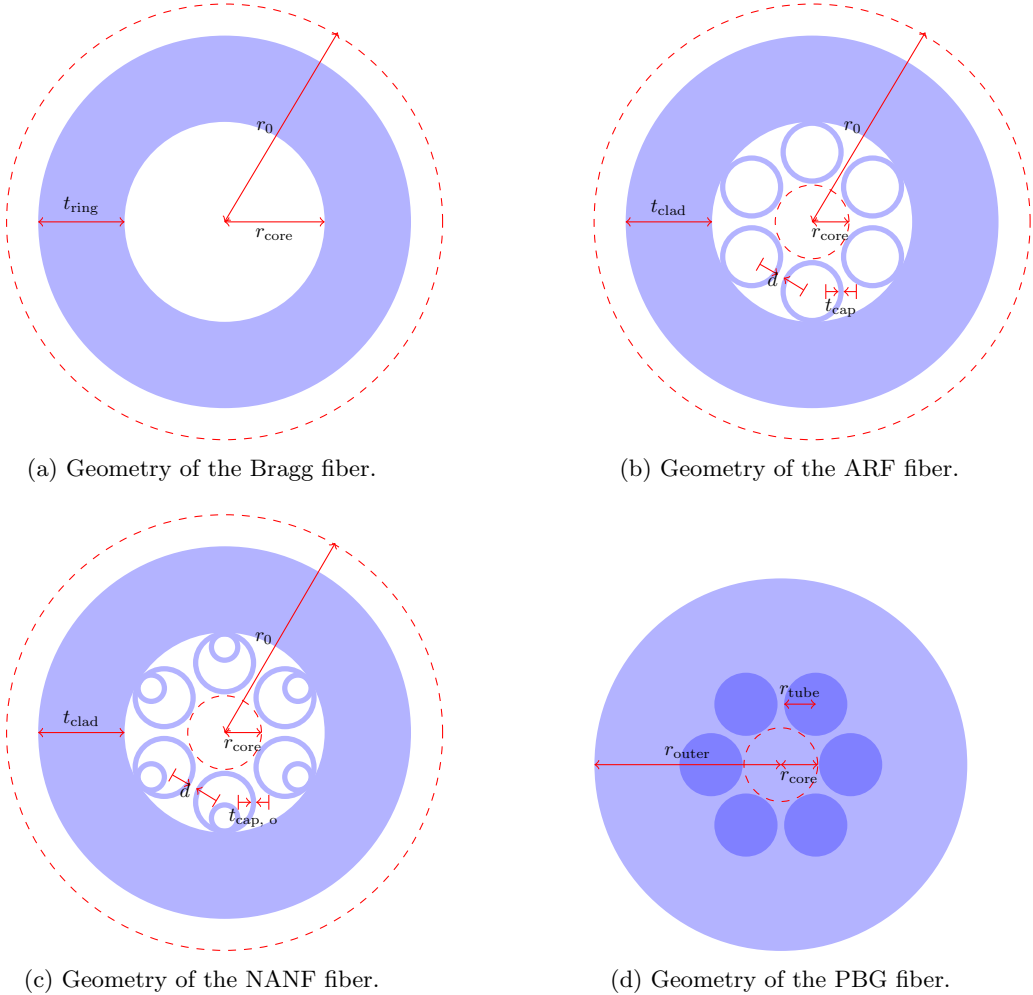

Modes are not available in analytic closed form for the three other fibers we consider. 
Their geometries are displayed in Figures~\ref{fig:geo_arf}, \ref{fig:geo_nanf} and~\ref{fig:geo_pbg}. 
Antiresonant fibers (ARF) such as those displayed in Figures~\ref{fig:geo_arf} and \ref{fig:geo_nanf} hold great promise due to their potential for low loss propagation in their large hollow cores (filled with air which has negligible optical nonlinearities). 
Negative-curvature hollow-core ARF, like that in Figure~\ref{fig:geo_arf}, but with eight thin capillaries were studied experimentally in~\cite{2013KolKosPryBirPlotDia}. 
The Nested Anti-Resonant Nodeless Fiber (NANF), whose microstructure design is illustrated in Figure~\ref{fig:geo_nanf}, has capillaries each containing a smaller capillary ring within it, was studied in~\cite{2014Pol}, where it was projected to have improved (lower) losses compared to the previous ARF designs, like that in Figure~\ref{fig:geo_arf}. 
Photonic bandgap (PBG) fibers like that  in Figure~\ref{fig:geo_pbg} form another category of popular microstructure fiber designs~\cite{Knigh03,OhterHirosYamad14, 2013PolPetRic, 1999CreManKniBirRusRobAll, 2003LitDunUsnEggWhiMcPdeSte} gaining utility, particularly where all-solid fibers are preferred. 
The geometry of the PBG fiber used here contains a single layer of hexagonally arranged high-index inclusions, a design shown to have the same transmission minima as more complex multi-layer devices \cite{2003LitDunUsnEggWhiMcPdeSte}. 
Note that it is not the goal of this paper to discuss the relative merits of various microstructures. 
Rather, our purpose in considering the diverse designs in Figure~\ref{fig:geometries} is to demonstrate the general purpose utility of the proposed algorithm. 
As we shall see, the modes of all these fibers have fine-scale features that the algorithm is able to successfully resolve. 
 
The remainder of this paper is organized as follows. Section~\ref{sec:pmlmaxwell} describes the mathematical formulation of the eigenproblem starting from Maxwell equations. 
All details of the PML and the finite element discretization are included there. 
The DWR error estimator, its theoretical basis,  and the adaptivity algorithm appears in Section~\ref{sec:dwr-maxwell}. 
The methodology is verified in Section~\ref{sec:bragg} using a model Bragg fiber and its semi-analytically computable propagation constants. 
Finally, we conclude in Section~\ref{sec:three-fibers} by applying the algorithm to the above-mentioned three practically important fiber microstructures.

\section{The eigenproblem for leaky Maxwell modes}\label{sec:pmlmaxwell}

In this section, we quickly outline the steps to derive a weak formulation for outgoing leaky hybrid modes of an optical fiber using a perfectly matched layer. 
The resulting system, labeled~\eqref{eq:maxwell-weak-pml-0} below, is the eigenproblem we shall solve in the remaining sections. 
We explain how the variable coefficients in~\eqref{eq:maxwell-weak-pml-0} are derived since they are important for the implementation of the error estimators discussed later.

The governing equation for light propagation in an optical fiber is the time-harmonic Maxwell system for the electric field, \(\hat E\), and the magnetic field, \(\hat H\), where all time variations are of the form \(e^{-\imath \omega t}\) for some frequency $\omega > 0$, 
\begin{subequations}
	\label{eq:MaxwellFull}
	\begin{align}
		-\imath \omega \mu \hat H + \nabla \times \hat E & = 0, \\
		\imath \omega \varepsilon \hat E + \nabla \times \hat H & = 0,
	\end{align}
\end{subequations}
for all \( (x_0,x_1,x_2)\in \RR^3\), where $+x_2$ is the field's propagation direction, \(\varepsilon\) represents the electric permittivity and \(\mu\) denotes the magnetic permeability, and \(\nabla \times \cdot \) is the three-dimensional curl operator. 
We assume that \(\mu > 0\) is isotropic and constant, an assumption that holds well for most optical materials.  However,  \(\varepsilon\) may have anisotropies of the following form 
\(
	\varepsilon =  \mathrm{diag}( \veps_\tau, \veps_\tau, \veps_2) = 
		\begin{bsmallmatrix}
			\varepsilon_\tau \delta & 0 \\ 0 & \varepsilon_2
	\end{bsmallmatrix}
\)
where $\varepsilon_\tau$ and $\varepsilon_2$ are positive scalar functions and $\delta$ denotes the $2 \times 2$ identity matrix. 
The assumption of isotropy in $x_0x_1$~plane is made only for simplicity and can be easily relaxed for the considerations in this paper. 

We are interested in Maxwell solutions that propagate along the longitudinal direction $x_2$, which represents the direction of translational symmetry of our optical fiber. 
Thus, we seek solutions of the form
\begin{equation}
	\label{eq:ansatz}
	\hat E(x_0,x_1,x_2) = E(x_0,x_1) \, e^{\imath\beta x_2}, \qquad
	\hat H(x_0,x_1,x_2) = H(x_0,x_1) \, e^{\imath\beta x_2},
\end{equation}
for some {\em propagation constant} $\beta$ and vector fields $E, H$ to be determined. 
By using unit vectors \(e_0, e_1, e_2\in \RR^3\) in the coordinate directions, the {\em transverse} part of a vector field \(F = F(x_0, x_1, x_2)\) is its projection into the span of $e_0$ and $e_1,$ denoted by $F_\tau$, so that $F = F_\tau + F_2 e_2$ where $F_2$ is the {\em longitudinal} component of $F$. 
Let $R$ denote the 90~degree rotation operator that maps any transverse field $F_\tau = F_0 e_0 + F_1 e_1$ to $R F_\tau = F_1 e_0 - F_0 e_1$. 
Define the divergence and the scalar-valued curl for transverse fields by $\dive F_\tau = \partial_0 F_0 + \partial_1 F_1$ and \(\curl(F_\tau) = \diver(R F_\tau).\) The rotated gradient of a scalar field $\phi$ is defined by \(\rot \phi = R \ \grad(\phi)\). 
Next, we substitute~\eqref{eq:ansatz} into~\eqref{eq:MaxwellFull}, decompose $E$ and $H$ into transverse and longitudinal components, and use the identity $\nabla \times \left(u\, e^{\imath\beta x_2}\right) = \left( e^{\imath\beta x_2} \curl u_\tau\right)e_2 + e^{\imath\beta x_2}\left( \rot u_2 - \imath\beta R u_\tau \right)$. 
Then, eliminating $H_2$ and $H_\tau$ and simplifying (see e.g., \cite[Appendix~A.1]{2023VanGopGro} for more details), the Maxwell system reduces to a system coupling the transverse electric field $E_\tau$ and a scaled longitudinal component of the field \(\varphi = \imath \beta E_2\), which reads as follows: 
\begin{subequations}
	\label{eq:maxwell-nonnondim-2}
	\begin{align}
		\label{eq:6a}
		\rot\curl E_\tau - \omega^2\veps_\tau \mu E_\tau + \grad \varphi & = - \beta^2 E_\tau, \\
		\label{eq:6b}
		\veps_2 \mu \varphi + \dive (\veps_\tau \mu E_\tau) & = 0,
	\end{align}
\end{subequations}
This is a mixed eigensystem for the eigenvalue $-\beta^2$ where the second equation may be viewed as a constraint. 

The physical dimensions of the fiber's transverse cross-section are in micrometers, while practical propagation constants are usually of the order of $10^6$, so a nondimensionalization is useful. 
Let \(\varepsilon_0\) and \(\mu_0\) be the electric permittivity and magnetic permeability of the vacuum, respectively. 
Let \(k^2 = \omega^2 \varepsilon_0 \mu_0\). 
Define the transverse and longitudinal refractive indices by \(n_\tau^2 = \varepsilon_\tau \mu_0 / \varepsilon_0 \mu \) and \(n_2^2 = \varepsilon_2 \mu_0 / \varepsilon_0 \mu\). 
Now, using a characteristic length scale \(L\) for the fiber cross-section, we rescale \((x_0, x_1) \mapsto \left(x_0/L, x_1/L\right)\) and consider all the
unknown functions in the rescaled coordinates. 
Outside some finite radius (of the same order as $L$) in the transverse plane we assume that the refractive index $n(x_0, x_1)$ is isotropic and equal to a constant
$n_0$. 
Using an {\em index well} defined by $V(x_0, x_1) = L^2 k^2 (n_0^2  - n_\tau(x_0, x_1)^2)$ in the rescaled nondimensional coordinates, the system~\eqref{eq:maxwell-nonnondim-2} leads us~\cite[Appendix~A.2]{2023VanGopGro} to the following non-dimensional equations: 
\begin{subequations}
	\label{eq:maxwell-nondim-0}
	\begin{align}
 		\rot \curl E_\tau + V  E_\tau + \grad  \varphi & = Z^2 E_\tau,
		&& \text { in } \RR^2,
		\\
		n_2^2  \varphi + \diver( n_\tau^2  E_\tau ) & = 0,
		&& \text { in } \RR^2,
	\end{align}
\end{subequations}
where $Z^2 = L^2 (k^2 n_0^2 - \beta^2)$ is a nondimensional eigenvalue. 

The system~\eqref{eq:maxwell-nondim-0} must be supplemented with boundary conditions. \rev{We shall use the outgoing condition at infinity, described below. 
Before discussing it however,}  consider the case where all fields have decayed to some negligible magnitude (as in a guided mode) at the boundary of a large enough disk $\om$ in the $x_0x_1$-plane, sufficiently far from the compact support of the index well~$V$. 
Multiplying the first equation of~\eqref{eq:maxwell-nondim-0} by a vector test function $F$ and the second by a scalar test function $\psi$, integrating by parts, and using $E_\tau =0$ and $ \varphi = 0$ on $\partial \om$, we obtain the following weak formulation. 
Find $E_\tau \in \Ho(\curl, \om)$ and $\varphi \in \Ho^1(\om)$, satisfying 
\begin{subequations}
	\label{eq:weak}
	\begin{align}
		\label{eq:weak-a}
          \int_\om(\curl E_\tau)( \curl F)  + \int_\om
          V E_\tau \cdot  F +  \int_\om(\grad \varphi) \cdot F
          & = Z^2 \int_\om E_\tau \cdot  F, 
          \\
          \label{eq:weak-b}
          \int_\om n_2^2 \varphi \psi - \int_\om n_\tau^2 E_\tau \cdot \grad \psi
          & = 0,
	\end{align}  
\end{subequations}
for all $F \in \Ho(\curl, \om)$ and all $\psi \in \Ho^1(\om)$. 
Here $\Ho(\curl, \om)$ denotes the space of square-integrable vector fields whose $\curl$ is also square integrable and whose tangential component along $\partial\om$ vanishes and $\Ho^1(\om)$ denotes square-integrable scalar functions which vanish on $\partial \om$ all of whose first-order derivatives are also square integrable on $\om$. 
Prior works featuring similar weak formulations include~\cite{KoshiInoue92, LeeSunCende91, 2023VanGopGro, VardaDemko03}. 

Leaky modes satisfy~\eqref{eq:maxwell-nondim-0} and also satisfy the
condition at infinity that the mode must be {\em outgoing}.  \rev{To
  explain this, first note that the fiber's translational symmetry in
  $x_2$-direction and the absence of inhomogeneities in
  $r^2\equiv x_0^2 + x_1^2 > R^2$ for some finite cylindrical radius
  $R$, implies that we can use separation of variables to obtain a
  series expansion of solution components in terms of solutions of the
  Bessel equation in $r$ and sines and cosines in the angular variable
  $\theta$.  From among all solutions to the Bessel equation, when $Z$
  is a positive real number, the outgoing requirement selects Hankel
  functions of the first kind, $H_\ell^{(1)}$. Hence the radial
  variation of all terms in the solution
  representation are of the form
  $H_\ell^{(1)}(Z r)$ as $r \to \infty$ (see e.g.,~\cite{ColliMonk98} or \cite[eq.~(23)]{2022GopParVan}).  For complex values of $Z$,
  the outgoing condition  requires
  that the solution is expressible in terms of the analytic
  continuations into $\mathbb C$, from the positive real axis,
  of the same Hankel functions $H_\ell^{(1)}(Z r)$.}
  A practically convenient method to impose
this condition is through a perfectly matched layer (PML)
\cite{1994Ber} which transforms coordinates so that the solution
of~\eqref{eq:maxwell-nondim-0}, while unaltered in a bounded region
where $V$ is inhomogeneous, becomes exponentially decaying outside.
The exponential decay allows us to truncate the infinite domain
$\RR^2$ in~\eqref{eq:maxwell-nondim-0} to a bounded computational
domain $\om$ as in~\eqref{eq:weak}.  To view PML as a complex
coordinate transformation~\cite{2003Mon, ColliMonk98}, we let
$\tilde x := \Phi(x)$ where $\Phi(x) = x\,\eta(r) / r$ for some
complex-valued function $\eta(r)$ to be specified later.  Here
\( x = [ x_0, x_1]^\mathsf{T} \),
\( \tilde x = [ \tilde x_0, \tilde x_1 ]^\mathsf{T} \), and
\(r=\sqrt{x_0^2 + x_1^2}\).  Put
\( x_\bot = [ x_1, -x_0]^\mathsf{T} \).  Consider the $2 \times 2$
Jacobian matrix whose entries are
$J_{ij} := \partial \tilde x_i / \partial x_j$.  Simple computations
give that with $\zeta(r) = \eta(r)/r$, we have
\begin{align}
	\label{eq:Jacobian}
	J & = \frac{{\zeta}'}{r} x x^\mathsf{T} + {\zeta} \delta, \qquad 
		J^{-1} = \frac{1}{\eta'} \delta + \frac{{\zeta}'}{\eta \eta'}x_\bot x_\bot^\mathsf{T}, \qquad
		\det(J) = {\zeta} \eta'.
\end{align}
We transform scalar-valued functions \(\psi\) and vector-valued functions \(\Psi\) in the following distinct ways:
\begin{equation}
	\label{eq:trafo-scalar-vector}
	\tilde \psi = \psi \circ \Phi^{-1}, \qquad
	\tilde \Psi = J^{-\mathsf{T}} \Psi \circ \Phi^{-1}.
\end{equation}
By placing a tilde over a differential operator, we indicate that the derivative is being taken with respect to $\tilde x$ rather than $x$. 
Calculations using the chain rule shows  that the matrix $ D \Psi$ whose $(i,j)$th entry equals $\partial_j \Psi_i $ satisfies $\tilde \partial_j \tilde \Psi_i = [ J^{-\mathsf{T}} \ D\Psi \ J^{-1} ]_{ij} +
	\sum_l \Psi_l\, \partial^2 x_l / \partial \tilde x_j \partial \tilde x_i$. 
Relating its skew symmetric part to the curl and applying~\eqref{eq:Jacobian}, we obtain 
\begin{subequations}
	\label{eq:curl-grad-map}
	\begin{equation}
		\tilde \curl \tilde \Psi = \frac{1}{(\eta')^2} \left( 1 + \frac{{\zeta}' r^2}{\eta} \right) (\curl \Psi ) \circ \Phi^{-1}
	\end{equation}
	Similarly, the chain rule applied to $\tilde \psi = \psi \circ \Phi^{-1}$ gives
	\begin{equation}
		\tilde\grad \,\tilde \psi = J^{-\mathsf{T}} (\grad \psi) \circ \Phi^{-1}.
	\end{equation}
\end{subequations}
In the complexified domain $\tilde \om = \Phi( \om)$, we solve the analogue of system~\eqref{eq:weak} obtained by replacing $\curl$ by $\tilde \curl$, $\grad$ by $\tilde \grad$, and $\om$ by $\tilde \om$, namely 
\begin{subequations}
	\label{eq:weak-2}
	\begin{align}
		\label{eq:weak-a-2}
          \int_{\tilde\om} (\tilde\curl \tilde E_\tau) ( \tilde\curl \tilde F)
          + \int_{\tilde\om} V \tilde E_\tau \cdot \tilde F
          + \int_{\tilde\om} (\grad \tilde\varphi ) \cdot  \tilde F)
          & = Z^2 \int_{\tilde\om} \tilde E_\tau \cdot  \tilde F,
          \\
          \label{eq:weak-b-2}
          \int_{\tilde\om} n_2^2 \tilde\varphi  \tilde\psi
          -
          \int_{\tilde\om} n_\tau^2 \tilde E_\tau\cdot \grad \tilde\psi
          & = 0.
	\end{align}  
\end{subequations}
where $\tilde E_\tau, \tilde F, \tilde \varphi$ and $ \tilde \psi$ are obtained from $E_\tau, F, \varphi$ and $ \psi$ using~\eqref{eq:trafo-scalar-vector}. 

Since it is more convenient to compute on a domain with real coordinates, we now transform the equations of~\eqref{eq:weak-2} back to $\om$ using the identities in~\eqref{eq:Jacobian} and \eqref{eq:curl-grad-map}, e.g., 
\[
  \int_{\tilde\om} (\tilde \curl \tilde E_\tau)  ( \tilde \curl \tilde F)
  =
  \int_{\tilde\om}  \kappa \,\curl E_\tau \, \curl F,
\]
with 
\begin{equation}
	\label{eq:kappa-defn}
	\kappa = \frac{{\zeta}}{(\eta')^3} \left( 1 + \frac{{\zeta}' r^2}{\eta} \right)^2.  
\end{equation}
Also defining 
\begin{equation}
	\label{eq:gamma-defn}
	\gamma = {\zeta} \eta' J^{-\mathsf{T}}J^{-1},
\end{equation}
and performing similar transformations for other terms in~\eqref{eq:weak-2}, we are led to a weak formulation with PML. 
Namely, we find $E_\tau \in \Ho(\curl, \om)$ and $\varphi \in \Ho^1(\om)$, satisfying 
\begin{subequations}
	\label{eq:maxwell-weak-pml-0}
	\begin{align}
		a(E_\tau,F) + c(\varphi,F) &    = Z^2\ m(E_\tau,F),
		&&\text{ for all } F \in \Ho(\curl,\om),
		\\
		b(E_\tau,\psi) + d(\varphi,\psi) &    = 0,
		&&\text{ for all } \psi \in \Ho^1(\om),
	\end{align}
\end{subequations}
where 
\begin{subequations} \label{eq:abcdm}
	\begin{align} 
		a(E_\tau,F) & = \left( \kappa \curl E_\tau,  \curl F\right)_\om
		+
		(V \gamma E_\tau, F)_\om,
		\\
		b(E_\tau,\psi) & = \left(n_\tau^2 \gamma E_\tau,  \grad \psi\right)_\om,
		\qquad  
		c(\phi,F) = \left( \gamma \grad \phi, F\right)_\om, \qquad 
		\\
		d(\phi,\psi) & = - \left( n_2^2 {\zeta} \eta' \phi,\psi\right)_\om,
		\qquad \qquad 
		m(E_\tau,F) = \left( \gamma E_\tau,  F\right)_\om.
	\end{align}
\end{subequations}
Here $(\cdot, \cdot)_\om$ denotes the (complex) inner product of $L^2(\om)$ or its Cartesian products and we have used the conjugates of test functions in~\eqref{eq:weak-2}.
The system~\eqref{eq:maxwell-weak-pml-0} completes the description of the weak formulation for computing leaky modes, except for a prescription of $\eta(r)$.

There are multiple ways to choose  $\eta(r)$, as can be seen from the literature~\cite{ColliMonk98, 2022GopParVan, KimPasci09, NanneWess18}. 
We use a two-dimensional analogue of an expression in~\cite{KimPasci09}. 
To describe it, we fix $\om$ to be a disk of radius $r_1$, assume that support of $V$ is contained in a disk of radius $r_0 < r_1$, and that a cylindrical PML is set in the annular region \(r_0 < r < r_1\). 
Let \(0 < \alpha\) be the PML strength parameter. 
We set 
\begin{subequations}\label{eq:eta-defn}
	\begin{equation}
		\eta(r)  = 1 + \imath \phi(r),
	\end{equation}
	where 
	\begin{equation}
		\phi(r) =
			\begin{cases}
				0, & \text{ if } r < r_0, \\
				\alpha\, \displaystyle{
					\frac{\displaystyle{\int_{r_0}^r (s - r_0)^{2} (s - r_1)^{2} \,
					\mathrm{d}s}}{
					\displaystyle{
					\int_{r_0}^{r_1} (s - r_0)^{2} (s - r_1)^{2} \, \mathrm{d}s}
					},
				}
				& \text{ if } r_0 < r < r_1.
			\end{cases}
	\end{equation}
\end{subequations}
This together with  $\zeta(r) = \eta(r) / r$ defines all quantities that appear  in~\eqref{eq:maxwell-weak-pml-0}. 

We conclude this section by describing the finite element discretization of~\eqref{eq:maxwell-weak-pml-0} that we shall employ. 
First, mesh $\om$ by a geometrically conforming finite element mesh of triangles, and denote the mesh by $\oh$. 
On a triangle $K$, let $P_p(K)$ denote the space of polynomials in two variables of degree at most~$p$. 
The degree $p$ and the maximal element diameter $h = \max_{K \in \oh} \mathrm{diam}(K)$ determine the richness of the discretization. 
The Lagrange finite element space on $\oh$, for any $p\ge 0$ is defined by $V_h = \{ \psi: \psi|_K \in P_{p+1}(K)$ for all $K \in \oh,$ $\psi$ is continuous$\}$. 
The N\'{e}d\'{e}lec finite element space~\cite{2003Mon} is the set of all vector fields $F$ on $\om$ whose tangential components are continuous across element interfaces and whose restrictions to each element $K \in \oh $ is a polynomial of the form $F|_K \in P_p(K)^2 + \left[\begin{smallmatrix} \phantom{-}x_1\\ -x_0 \end{smallmatrix}\right] P_p(K)$. 
We use curved triangles next to curved material interfaces, in which case, as usual, the polynomial space within such an element is replaced by the pullback of the above-indicated polynomial spaces from a reference triangle. 
The Galerkin discretization using these spaces give the discrete eigenproblem that solves for $\Eth \in N_h$ and $\varphi_h \in V_h$ such that 
\begin{subequations}
	\label{eq:maxwell-weak-pml-h}
	\begin{align}
		a(\Eth,F) + c(\varphi_h,F) &    = Z_h^2\ m(\Eth,F),
		&&\text{ for all } F \in N_h, 
		\\
		b(\Eth,\psi) + d(\varphi_h,\psi) &    = 0,
		&&\text{ for all } \psi \in V_h.
	\end{align}
\end{subequations}
All results of computations described in later sections are obtained by solving~\eqref{eq:maxwell-weak-pml-h}.

\section{The error estimator for eigenvalues}\label{sec:dwr-maxwell}

In this section, we describe the {\it a posteriori} error estimator for the previously described eigenproblem. 
It is obtained by applying the DWR technique~\cite{2001BecRan, 2001HeuRan, RannaWesteWolln10}. 
Although only eigenproblems with second-order elliptic operators were considered in these references, the generality of their technique is well recognized. 
Even though our eigenproblem involves a non-elliptic operator, the same approach applies, as we shall now see. 

Let $H = \Ho(\curl, \om) \times \Ho^1(\om)$ and $ H_h = N_h \times V_h.$ Putting $\lambda=Z^2$, $u = (E_\tau, \varphi) \in H$, and $v = (F, \psi) \in H$, we restate
the eigenproblem~\eqref{eq:maxwell-weak-pml-0} as 
\begin{align} \label{eq:gen_eig_weak-0}
  A(u,w)
  & = \lambda B(u,w),
    \qquad 
  \| u \|_H  = 1,
\end{align}
for all $w \in H$ with 
\begin{gather*}
  A(u, v) =
  a(E_\tau,F) + c(\varphi,F)
  + b(E_\tau,\psi) + d(\varphi,\psi),
  \quad
  B(u, v) = m(E_\tau, F),
\end{gather*}
with $a, b, c, d, m$ as in \eqref{eq:abcdm}. 
Viewing $u$ as a right eigenfunction, we also assume that there is a corresponding left eigenfunction $0 \ne \ut = (\tilde{E}_\tau, \tilde \varphi) \in H$ satisfying 
\begin{subequations} \label{eq:gen_eig_adj-weak-0}
  \begin{align}
    \label{eq:gen_eig_adj-weak-0-eq}
    A(\vt, \ut) & =  \lambda B( \vt, \ut), \quad \text{ for all } \vt \in H,\\
    \label{eq:gen_eig_adj-weak-0-normalize}
    B(u, \ut) & = 1.
  \end{align}
\end{subequations}
Next, consider the Galerkin discretizations of these eigenproblems, considered previously in~\eqref{eq:maxwell-weak-pml-h}. 
Using the same type of normalizations as in~\eqref{eq:gen_eig_weak-0}--\eqref{eq:gen_eig_adj-weak-0}, we assume there is a discrete right eigenfunction $u_h = (\Eth, \varphi_h) \in H_h$ and a discrete left eigenfunction $\ut_h = (\tEth, \tilde\varphi_h) \in H_h$ satisfying 
\begin{align} \label{eq:gen_eig_weak-h}
  A(u_h,w)
  & = \lambda_h B(u_h,w), \qquad
    \| u_h \|_H = 1,
\end{align}  
for all $w \in H_h$ and 
\begin{subequations} \label{eq:gen_eig_adj-weak-h}
  \begin{align}
    \label{eq:gen_eig_adj-weak-h-eq}
    A(\vt, \ut_h)
    & =  \lambda_h B( \vt, \ut_h), \quad \text{ for all } \vt \in H_h,\\
    \label{eq:gen_eig_adj-weak-h-normalize}
    B(u_h, \ut_h)
    & = 1.
  \end{align}
\end{subequations}
where $\lambda_h = Z_h^2$. 

Next, to describe the error estimators, using the coefficients $\kappa$ and $\gamma$ in~\eqref{eq:kappa-defn} and \eqref{eq:gamma-defn}, we first set element-wise residuals. 
Let $\| \cdot \|_D$ denote the $L^2(D)$-norm for a subset $D$ such as $T \in \oh$ or its boundary $\partial T$. 
Let $h_T$ denote the diameter of $T$.  
On an interior mesh edge, let $\nu$ denote a unit normal of arbitrarily fixed orientation. 
Jumps of multivalued functions at the element interfaces are denoted by $\llbracket \cdot \rrbracket$; so, for example, $\llbracket (n_\tau^2\gamma \Eth) \cdot \nu \rrbracket$ denotes the jump of the normal ($\nu$) component of $n_\tau^2\gamma \Eth$ across an interior edge. 
Using this notation, define (unweighted) element-wise residuals by 
\begin{subequations}
	\label{eq:ee-eta}
\begin{align}
    \rho_{T, 1}^2 
    & =
      h_T^2\left\lVert
      \rot (\kappa  \curl \Eth)  + 
      V\gamma  \Eth +
      \gamma\grad \varphi_h -
      Z_h^2\gamma \Eth
      \right\rVert_T^2
  \\  \nonumber 
  & \quad +
      \frac{h_T}{2} 
      \big\lVert
      \left\llbracket
      \kappa \curl \Eth  \right\rrbracket
      \big\rVert_{\partial T \setminus \partial \om}^2,
  \\
  \rho_{T, 2}^2
    & = h_T^2 \left\| \diver( \gamma \grad \varphi_h + V \gamma \Eth 
      - Z_h^2 \gamma \Eth) \right\|_T^2
  \\ \nonumber 
  & \quad + \frac{h_T}{2} \left\| \left\llbracket
      (  \gamma \grad \varphi_h + V \gamma \Eth 
      - Z_h^2 \gamma \Eth) \cdot \nu
    \right\rrbracket \right\|_{\partial T \setminus \partial \om}^2,
  \\
  \rho_{T, 3}^2
    & = h_T^2
      \left\| \diver( n_\tau^2 \gamma \Eth) + n_2^2 \zeta \eta' \varphi_h
      \right\|_T^2
      +
      \frac{h_T}{2}
      \left\| \llbracket
      (n_\tau^2 \gamma \Eth) \cdot \nu
      \rrbracket
      \right\|_{\partial T \setminus \partial\om}^2,
\end{align}
and define weights arising from the dual eigenfunction components by 
\begin{align}
  \tilde\omega_{T, 1}    =  \| \curl \tilde E_\tau\|_T, \qquad 
  \tilde\omega_{T, 2}   =  \| \tilde E_\tau\|_T, \qquad
  \tilde\omega_{T, 3}  = \| \grad \tilde\varphi \|_T.  
\end{align}
Similarly, using complex conjugates of the coefficients, we define analogous element-wise residuals for the dual problem, 
\begin{align}
  \tilde\rho_{T, 1}^2
  & =
    h_T^2\left\lVert
    \rot (\bar{\kappa}  \curl \tEth)  +
    n_\tau^2 \bar\gamma \grad \tilde\varphi_h + 
    V\bar\gamma  \tEth
    - \bar Z_h^2 \bar \gamma \tEth
    \right\rVert_T^2
  \\ \nonumber 
  & \quad +
    \frac{h_T}{2} 
    \Big\lVert
    \left\llbracket
    \bar \kappa \curl \tEth  \right\rrbracket
    \Big\rVert_{\partial T \setminus \partial \om}^2,
  \\
  \tilde\rho_{T, 2}^2
  & =
    h_T^2 \left\| \diver(n_\tau^2\bar\gamma \grad \tilde\varphi_h
    + V \bar \gamma \tEth
    - \bar Z_h^2  \bar \gamma \tEth)
    \right\|_T^2
  \\ \nonumber 
  & \quad +
    \frac{h_T}{2} 
    \Big\lVert
    \left\llbracket
    (n_\tau^2\bar\gamma \grad \tilde\varphi_h
    + V \bar \gamma \tEth
    - \bar Z_h^2  \bar \gamma \tEth
    )\cdot \nu
    \right\rrbracket
    \Big\rVert_{\partial T \setminus \partial \om}^2,
  \\
  \tilde\rho_{T, 3}^2
  & =
    h_T^2 \left\| \diver( \bar \gamma \tEth)
    + n_2^2 \bar\zeta \bar\eta' \tilde\varphi_h
    \right\|_T^2
    +
    \frac{h_T}{2}
    \Big\lVert
    \left\llbracket
    (\bar\gamma \tEth ) \cdot \nu
    \right\rrbracket
    \Big\rVert_{\partial T \setminus \partial \om}^2,    
\end{align}
and weights generated by the right eigenfunction: 
\begin{align}
  \omega_{T, 1}    = \| \curl  E_\tau\|_T, \qquad 
  \omega_{T, 2}   =   \| E_\tau\|_T, \qquad
  \omega_{T, 3}  = \| \grad \varphi \|_T.  
\end{align}
Using these, we define a local element-wise error indicator by 
\begin{equation}
  \label{eq:eta}
  \eta_T =
  \sum_{i=1}^3
  \rho_{T, i} \,\tilde\omega_{T, i}
  +
  \tilde\rho_{T, i} \,\omega_{T, i},  
\end{equation}
\end{subequations}
for each $T \in \oh$. A practical version of this idealized error indicator will be given later: see~\eqref{eq:eta-h}.

The adaptive algorithm we use is motivated by the next result (Theorem~\ref{thm:drw_error-maxwell} below), which gives a global reliability estimate for the error indicator $\eta_T$. 
We prove it using the DWR approach laid out in~\cite{2001BecRan, 2001HeuRan, RannaWesteWolln10}, but bring in an additional ingredient to decompose $H(\curl)$ approximation errors into locally controllable gradients and a regular remainder using the results of~\cite{Schob08b}. 
As usual, we make a ``saturation assumption,'' which here takes the form that the discretizations considered are fine enough so that 
\begin{equation}
  \label{eq:saturation}
  \rev{
  \left(\max_{x \in \om} \| \gamma (x)\|_2\right)
  \,\| E_\tau - \Eth \|_{\om}\;
    \| \tilde E_\tau - \tEth \|_{\om}\,
    < 1
    }
\end{equation}
holds. 
The assumption is likely to  hold on reasonable  meshes since the product of errors on the left is expected to go to zero faster than the field error as
$h$ approaches zero. 
Hereon, we use $A \lesssim B $ to indicate that there exists a meshsize ($h$) independent constant $C>0$ (possibly dependent on shape regularity of the mesh) such that the inequality $A \le C B$ holds. 

\begin{theorem}
  \label{thm:drw_error-maxwell}
  Suppose the vectors $u=(E_\tau, \varphi)$, $u_h=(\Eth, \varphi_h),$
  $\tilde u=(\tilde E_\tau, \tilde\varphi)$,
  $\tilde u_h=(\tEth, \tilde\varphi_h)$,
  solve~\eqref{eq:gen_eig_weak-0}, \eqref{eq:gen_eig_adj-weak-0},
  \eqref{eq:gen_eig_weak-h} and \eqref{eq:gen_eig_adj-weak-h},
  respectively, with accompanying exact and discrete eigenvalues
  $\lambda$ and $\lambda_h$. Suppose also that~\eqref{eq:saturation}
  holds. Then the error in eigenvalue can be bounded by the error indicators defined in~\eqref{eq:ee-eta} by 
  \begin{equation}
    \abs{\lambda - \lambda_h} \,\lesssim\,
    \sum_{T \in \Omega_h} \eta_T.
  \end{equation}
\end{theorem}
\begin{proof}
  Let
  \[
    \rho(\vt) = A(u_h, \vt) - \lambda_h B(u_h, \vt), \qquad 
    \rhot(v) = A(v, \ut_h) - \lambda_h B(v, \ut_h)
  \]
  for any $v ,\vt \in H$. Also let
  $\sigma = B (u - u_h, \ut - \ut_h)/2$.  The argument is
  based~\cite{2001HeuRan} on the identity
  \begin{equation}
    \label{eq:key-identity}
    (\lambda - \lambda_h) (1 - \sigma) = 
    \frac 1 2\rho( \ut- \vt_h) +
    \frac 1 2\rhot( u - v_h)
  \end{equation}
  which holds for any $\vt_h, v_h \in H_h$. Its proof is elementary:
  \begin{align*}
    \rho( \ut- \vt_h) +
    \rhot( u - v_h)
    & = A(u_h, \tilde u) - \lambda_h B(u_h, \ut)
    \\
    & + A(u , \ut_h) - \lambda_h B(u, \ut_h),
    && \text{ by \eqref{eq:gen_eig_weak-h} \& \eqref{eq:gen_eig_adj-weak-h},}
    \\
    & = (\lambda - \lambda_h)
      \big[ B(u_h, \ut) + B(u, \ut_h) \big],
    && \text{ by \eqref{eq:gen_eig_weak-0} \& \eqref{eq:gen_eig_adj-weak-0},}
    \\
    & = (\lambda - \lambda_h)
      \big[ B(u_h -u, \ut - \ut_h) + 2\big],
    && \text{ by \eqref{eq:gen_eig_adj-weak-0-normalize} \& \eqref{eq:gen_eig_adj-weak-h-normalize},}
  \end{align*}
  from which~\eqref{eq:key-identity} follows for any $v_h, \tilde v_h \in H_h$.

  We proceed setting
  $v_h =(\varPi_h E_\tau, I_h \varphi_h)$, where $I_h$ is the
  Scott-Zhang interpolant \cite{ScottZhang90} and $\varPi_h$ is the
  quasi-interpolant defined in~\cite{Schob08b}. The results
  of~\cite{Schob08b} prove  that there is an $e \in \Ho^1(\om)$ and
  $z \in [\Ho^1(\om)]^2$ such that
  \begin{subequations}
    \label{eq:E-decom-e-z}
    \begin{align}
      \label{eq:E-decom}
      E_\tau - \varPi_h E_\tau & = z + \grad e,  \\
      \label{eq:e-bd}
      h_T^{-1} \| e \|_T + \| \grad e \|_T & \lesssim \| E \|_{\om_T},\\
      \label{eq:z-bd}
      h_T^{-1} \| z \|_T + \| \grad z \|_T & \lesssim \| \curl E \|_{\om_T}.
    \end{align}
  \end{subequations}
  where $\om_T$ is the union of all elements connected to $T$.
  Using this decomposition and  integrating by parts,
  \begin{align}
    \nonumber 
    \tilde\rho(u - v_h)
    & = A(u - v_h, \tilde u_h) - \lambda_h B( u - v_h, \tilde u_h)
    \\ \nonumber 
    & = A((z + \grad e, \varphi - I_h \varphi),  (\tEth, \tilde \varphi_h))
      -Z_h^2 \left(\gamma (z + \grad e), \tEth\right)_\om
    \\ \nonumber 
    & =
      \left(\kappa \curl z, \curl \tEth\right)_\om
    \\ \nonumber 
    & \qquad +
      \left(V \gamma (z + \grad e)
      +
      \gamma \grad( \varphi - I_h \varphi)
      - Z_h^2 \gamma (z + \grad e), \tEth\right)_\om
    \\ \nonumber 
    &\qquad + \left(n_\tau^2 \gamma (z + \grad e),
      \grad \tilde\varphi_h\right)_\om
      -\left(n_2^2 \zeta \eta' ( \varphi - I_h \varphi),
      \tilde \varphi_h\right)_\om
    %
    \\  \label{eg:sum-tilde-rho}
    & = \sum_{T \in \oh}
      \Big[
      \left(z \cdot t, \bar \kappa \curl \tEth\right)_{\partial T}
    \\ \nonumber 
    & \qquad\qquad +
      \left(z, \rot \bar\kappa \curl \tEth
      + n_\tau^2 \bar\gamma \grad \tilde\varphi_h + V \bar \gamma \tEth
      - \bar Z_h^2 \bar\gamma \tEth\right)_T      
    \\ \nonumber 
    &\qquad\qquad
      +
     \left(e,
      (n_\tau^2 \bar \gamma \grad \tilde\varphi_h + V \bar \gamma \tEth
      - \bar Z_h^2 \bar\gamma\tEth)\cdot \nu\right)_{\partial T}
    \\ \nonumber 
    &
      \qquad\qquad
      - \left(e, \diver(
      n_\tau^2 \bar \gamma \grad \tilde\varphi_h
      + V \bar \gamma \tEth 
      -\bar Z_h^2 \bar\gamma \tEth)\right)_T
    \\ \nonumber 
    & \qquad\qquad
      +(\varphi - I_h \varphi, (\bar \gamma \tEth)
      \cdot \nu )_{\partial T \setminus \partial \om}
    \\ \nonumber
    & \qquad\qquad
      - (\varphi - I_h \varphi,
      \diver(\bar \gamma \tEth) + n_2^2 \bar \zeta \bar \eta' \tilde \varphi_h)_T
      \Big],
  \end{align}
  where $t$ and $\nu$ are, respectively, the unit counterclockwise
  tangent and the unit outward normal vectors on element boundaries. To
  bound the element boundary terms in the sum above, we use the
  well-known local trace estimate
  $ h_T^{-1} \| w\|_{\partial T}^2 \lesssim h_T^{-2} \| w \|_T^2 + \| \grad w
  \|_T^2$ which holds for any $w \in H^1(T)$ on any element $T$
  of a shape-regular mesh. Using it, the first
  term in the sum~\eqref{eg:sum-tilde-rho} can be bounded by
  \begin{align*}
    \Big|  \sum_{T \in \oh}
    (z \cdot t,    \bar \kappa \curl \tEth)_{\partial T} \Big|
    &\lesssim   \sum_{T \in \oh}
    \left(h_T^{-1/2}\|z\|_{\partial T}\right) \,
      \left(h_T^{1/2}\left\| \jump{\bar \kappa \curl \tEth}
      \right\|_{\partial T \setminus \partial \om}\right)
    \\
    &\lesssim   \sum_{T \in \oh}
      (h_T^{-1} \| z \| + \| \grad z\|_T) \; \tilde \rho_{T, 1}
      \\
    &\lesssim 
      \sum_{T \in \oh}
      \| \curl E_\tau \|_{\om_T}\; \tilde \rho_{T, 1},
  \end{align*}
  where we have used~\eqref{eq:z-bd}.  The second term is also bounded using~\eqref{eq:z-bd}:
  \begin{align*}
    \big(z, \rot \bar\kappa \curl \tEth
    + n_\tau^2 \bar\gamma \grad \tilde\varphi_h + V \bar \gamma \tEth
    &
    - \bar Z_h^2 \bar\gamma \tEth\big)_T
      \;\lesssim \;
      \| \curl  E_\tau \|_{\om_T}\; \tilde \rho_{T, 1}.    
  \end{align*}
  Similarly, applying~\eqref{eq:e-bd}, we find
  \begin{align*}
    \Big|  \sum_{T \in \oh}
    & \left(e,
      (n_\tau^2 \bar \gamma \grad \tilde\varphi_h + V \bar \gamma \tEth
      - \bar Z_h^2 \bar\gamma \tEth)\cdot \nu\right)_{\partial T}
      \Big|
    \\
    & \lesssim
    \sum_{T \in \oh}
    \left( h_T^{-1/2} \| e \|_{\partial T}\right)
    \left( h_T^{1/2}
      \left\| \left\llbracket
      (n_\tau^2 \bar \gamma \grad \tilde\varphi_h + V \bar \gamma \tEth
      - \bar Z_h^2 \bar\gamma \tEth) \cdot \nu \right\rrbracket
      \right\|_{\partial T}\right)
    \\
    & \lesssim
      \sum_{T \in \oh} \left( h_T^{-1} \| e \|_T + \| \grad e\| \right)\,
      \tilde \rho_{T, 2}
      \;\lesssim\;
      \sum_{T \in \oh} \| E_\tau\|_{\om_T} \; \tilde \rho_{T, 2}
  \end{align*}
  and
  \begin{align*}
    \big(e, \diver(
      n_\tau^2 \bar \gamma \grad \tilde\varphi_h
    & + V \bar \gamma \tEth -\bar Z_h^2 \bar\gamma \tEth)\big)_T
      \le \;
      h_T^{-1} \| e\|_T \,\tilde\rho_{T, 2}
      \;\lesssim\; \| E_\tau\|_{\om_T} \; \tilde \rho_{T, 2}.
  \end{align*}
  For the last two terms in~\eqref{eg:sum-tilde-rho}, we use the
  well-known property of $I_h$,
  \[
    h_T^{-1} \| \varphi - I_h \varphi \|_T +  \| \grad(\varphi - I_h \varphi) \|_T
    \lesssim \| \grad \varphi\|_{\om_T}
  \]
  in place of~\eqref{eq:E-decom-e-z} and apply Cauchy-Schwarz
  inequality in a similar fashion as above. Gathering bounds on all terms,
  we arrive at 
  \begin{equation}
    \label{eq:rho-tilde-bound}
        \tilde\rho
    ( u - v_h)
    \lesssim
    \sum_{T \in \oh} \omega_{T, 1} \; \tilde \rho_{T, 1} + \omega_{T,
      2} \; \tilde \rho_{T, 2}
    + \omega_{T,    3} \; \tilde \rho_{T, 3}.
  \end{equation}

  Next, consider the other residual, $\rho( \tilde u - \tilde v_h)$. Setting 
  $\tilde v_h = (\varPi_h \tilde E_\tau, I_h \tilde \varphi_h)$ and
  applying the same type of arguments, we 
  prove that  
  \begin{equation}
    \label{eq:rho-bound}
    \rho
    ( u - \tilde v_h)
    \lesssim
    \sum_{T \in \oh} \tilde\omega_{T, 1} \;  \rho_{T, 1}
    + \tilde\omega_{T,      2} \; \rho_{T, 2}
    + \tilde \omega_{T,    3} \;  \rho_{T, 3}.    
  \end{equation}
  \rev{Finally, note that the definitions of $\sigma$ and $B$ imply 
  \begin{align*}
    2 \sigma
    & = B (u - u_h, \ut - \ut_h) = m( E_\tau - \Eth, \tilde{E}_\tau - \tEth)
    \\
    & = (\gamma   ( E_\tau - \Eth), \tilde{E}_\tau - \tEth)_\om
      < 1,
  \end{align*}
  where the last inequality is due to~\eqref{eq:saturation}.
  Thus $1 - \sigma > 1/2$. Using this, together with~\eqref{eq:rho-bound} and~\eqref{eq:rho-tilde-bound} 
within the identity~\eqref{eq:key-identity},  the stated estimate follows.}
\end{proof}

We conclude this section by describing the {\em adaptive algorithm} we implemented. 
First, to make the error indicator $\eta_T$ practical, following the heuristics of~\cite[Remark~10]{2001HeuRan}, we take the next step of replacing the unknown weights $\tilde \omega_{T, i}$ and $\omega_{T, i}$ by a computable analogue that is likely to be close in value, namely
\begin{align*}
	\tilde\omega^{(h)}_{T, 1}    =  \| \curl \tEth\|_T, \qquad 
	\tilde\omega_{T, 2}^{(h)}   =  \| \tEth\|_T, \qquad                        
	\tilde\omega_{T, 3}^{(h)}  = \| \grad \tilde\varphi_h \|_T, \\
	\omega_{T, 1}^{(h)}    =   \| \curl  \Eth\|_T, \qquad 
	\omega_{T, 2}^{(h)}   = \| \Eth \|_T, \qquad                       
	\omega_{T, 3}^{(h)}  = \| \grad \varphi_h \|_T.  
\end{align*}
Replacing  the weights in~\eqref{eq:eta} by these, we obtain the practical element-wise error indicator we use: 
\begin{equation}
	\label{eq:eta-h}
	\eta^{(h)}_T(\lambda_h, u_h, \tilde u_h) =
	\sum_{i=1}^3
	\rho_{T, i} \,\tilde\omega_{T, i}^{(h)}
	+
	\tilde\rho_{T, i} \,\omega_{T, i}^{(h)}.
\end{equation}
\rev{%
  The next practicality involves aggregating estimators when multiple
  eigenvalues $\lambda_h^{(\ell, k)}$
  are found clustered,  
  together with their left and
  right eigenvectors
  $\ut_h^{(\ell, k)}, u_h^{(\ell, k)}, \; k=1, \dots, K_{\ell}$, 
  in the \(\ell\)-th adaptive iteration.  Then, we simply
  compute the element-wise maximum of the error indicators for each
  eigenvalue. Namely, we compute
   \(\eta_T^{(\ell, k)} = \eta_T^{(h)}(\lambda_h^{(\ell, k)},
   u_h^{(\ell, k)}, \ut_h^{(\ell, k)})
   \)
   using~\eqref{eq:eta-h} and set
    \begin{equation}
        \label{eq:eta-h-max}
        \eta_T^{(\ell)} =
        \max_{k=1, \dots, K_{\ell}} \eta_T^{(\ell, k)}.
    \end{equation}
    An element \(T'\in \oh^{(\ell)}\) is marked for refinement if
    \begin{equation}
        \label{eq:mark}
        \eta_{T'}^{(\ell)} > \theta \max_{T \in \oh} \eta_T^{(\ell)},
    \end{equation}
    where \(\theta\) is an input refinement threshold parameter.
This is used in Algorithm~\ref{alg:adaptive-long}, which summarizes the adaptive strategy we implemented. }

\begin{algorithm}[tb] 
    \caption{\rev{Adaptive algorithm for leaky modes}} 
	\label{alg:adaptive-long}

    \SetKwData{CriteriaParam}{\(\theta\)}
	
	\KwIn{Initial mesh $\om_h^{(0)},$ refinement threshold  $\theta$, an eigenvalue search region $D \subseteq \mathbb{C}$, and the maximum permitted number of degrees of freedom \(N_{\mathrm{max}}\).
	}
	\KwResult{
      	Final adaptively refined mesh $\oh^{({\ell})}$, and a cluster of eigenvalues $\varLambda_h^{({\ell})}$ contained in $D$ obtained using \(\oh^{({\ell})}\).        
	}
	\vspace{1.0em}
	\For{${\ell} \gets 0, 1, 2, \ldots$}{
        SOLVE: Assemble the system in \eqref{eq:gen_eig_weak-h} using \(\oh^{({\ell})}\) and {solve} it using a sparse eigensolver to obtain a cluster of  eigenvalues $\lambda_h^{({\ell}, k)} = (Z_h^{({\ell}, k)})^2$ contained in $D$, and the corresponding right and left eigenvectors $u_h^{({\ell}, k)}$ and  $\ut_h^{({\ell}, k)}$ for $k=1, \dots, K_{\ell}$.\;

        ESTIMATE: {Compute} error indicators \(\eta_T^{({\ell})}\) using \eqref{eq:eta-h-max} for each element \(T \in \oh^{({\ell})}\).\;
        
        MARK 
        elements for refinement based on \eqref{eq:mark} and input $\theta$.\;  
		
        REFINE 
        marked elements in  \(\oh^{({\ell})}\)  (as well as surrounding elements as needed to obtain mesh conformity) to generate the next mesh \(\oh^{({\ell}+1)}\).\;
		
		\If{number of degrees of freedom on \(\oh^{({\ell}+1)} > N_{\mathrm{max}}\)}{
			break\;
		}
	}
	\KwRet{$\oh^{({\ell})}$ and $\varLambda_h^{({\ell})}$\;}
\end{algorithm}


\rev{We conclude this section by specifying the parameters we used in
  Algorithm~\ref{alg:adaptive-long} and ensuing computations.  We use
  NGSolve~\cite{Schobother} for N\'{e}d\'{e}lec and Lagrange finite
  elements and for assembly of the eigensystem. It supports curved
  triangles, with a specifiable order of polynomial curving, as well
  as high-order N\'{e}d\'{e}lec and Lagrange elements.  }
Algorithm~\ref{alg:adaptive-long} can use any sparse eigensolver.  In
our computations, we use the FEAST contour integral
eigensolver~\cite{2009Pol, GopalGrubiOvall20a} (whose specific
application to leaky modes is also detailed
in~\cite[Section~A.4]{2022GopParVan}), which takes as input a search
region enclosed by a simple closed contour (denoted by $D$ in the
algorithm).  We use circular contours of small enough radius centered
at
\( \hat \lambda_h^{(n)} = \frac{1}{K_n}\sum_{k=1}^{K_n}
\lambda_h^{(n,k)}\), i.e., these centers track changes within each
adaptive iteration and $D$ may therefore vary with iteration
number~$\ell$.  How aggressively the refinements are made is
determined by the parameter $\theta$.
For our computations, we used \rev{\(\theta = 0.75\)}, a value determined by trial and error to give us  a good balance between localization of the refinements and a steady increase in the number of degrees of freedom (d.o.f.s). 
So as to represent all eigenmodes of the eigenvalue cluster on the same mesh, the algorithm uses a marking strategy based on the maximum of the error indicators generated by each eigenmode. Other marking strategies are possible, as in~\cite{1994Ver, 2001BecRan, 2001HeuRan}. 
As a stopping criterion, we set the maximum number of degrees of freedom to \(N_{\text{max}} = 2 \times 10^6\) in our numerical experiments. 

\section{Verification by Bragg fiber}\label{sec:bragg}

It is well known \cite{YehYariv76,LitchAbeelHeadl02} that an air waveguide surrounded by a higher index cladding can have modes that propagate energy and are primarily confined in the air core region. 
This design is a type of hollow-core Bragg fiber such that its transverse cross-section consists of a glass ring in infinite air (see Figure~\ref{fig:geo_bragg}). 
Such fibers were more recently re-examined in~\cite{2017Bir,2023VanGopGro} where the dielectric ring was viewed as an ``anti-resonant'' layer. 
In particular, the results of~\cite{2023VanGopGro} demonstrated that certain fine-scale features in the higher index region must be resolved to obtain confinement losses with acceptable accuracy. 
There, these fine features were also definitively shown to be not an artifact of numerical approximations, nor of the PML, since they were captured both by semi-analytic methods and by numerical computations. 

These prior studies make the Bragg configuration an ideal example for verifying our adaptive strategy. 
Because we can semi-analytically compute the exact eigenvalue $\lambda = Z^2$ for this geometry, we are able to study the history of actual eigenvalue errors as the adaptive mesh refinement proceeds. 
Moreover, starting with an unbiased mesh that disregards the above-mentioned prior knowledge of the fine-scale features in the modes, we investigate if the (deliberately blinded) adaptivity process is able to sense these features automatically and guide refinement to capture them. 

\begin{figure}[htbp]
	\centering
	\begin{subfigure}[t]{\TWOCOLWIDTH}
		\centering
        \includegraphics[trim = {20cm 8.9cm 20cm 10cm}, clip, width = \ONECOLWIDTH]{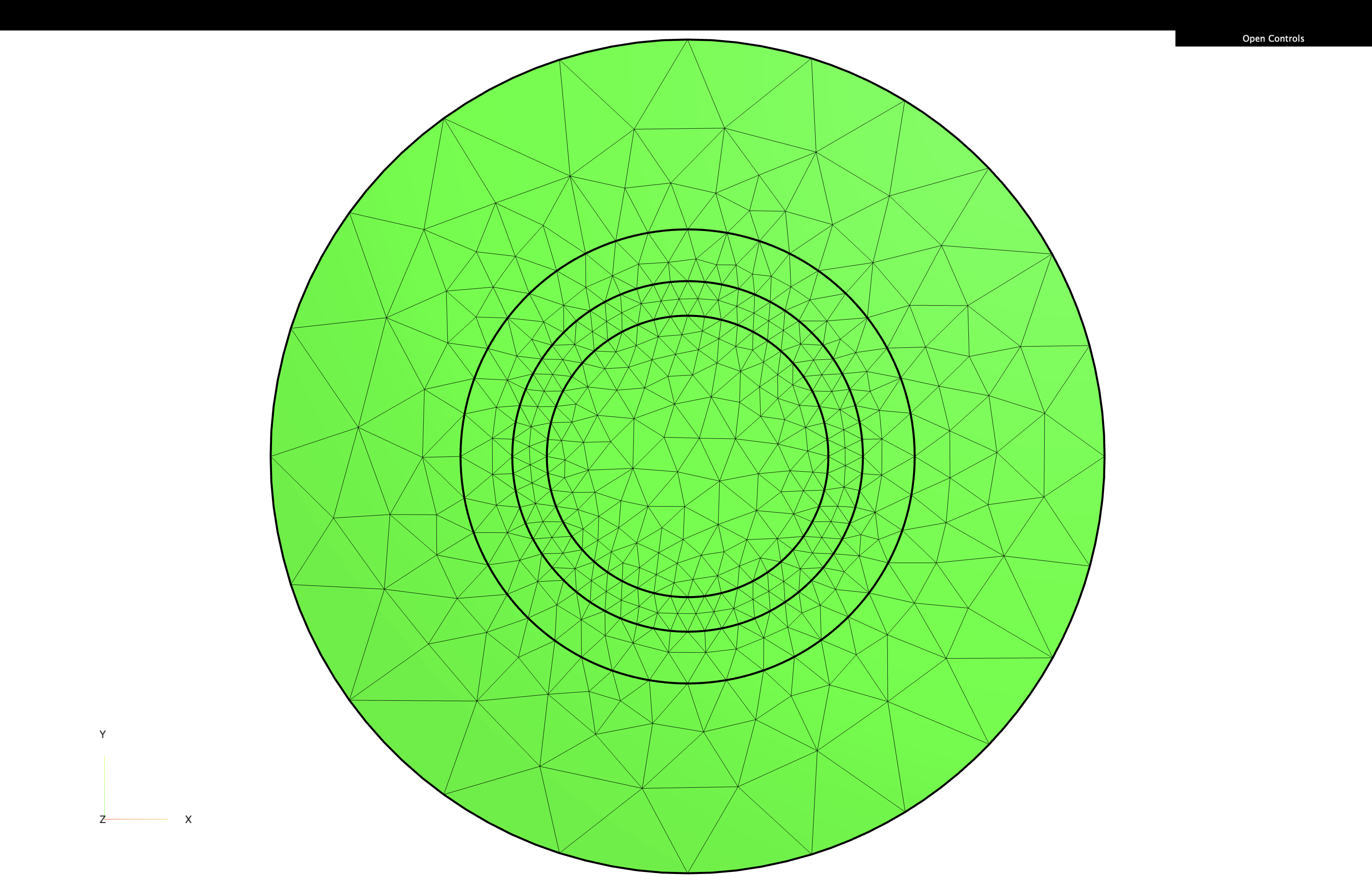}
		\caption{Initial mesh for the Bragg fiber, zoomed in near the glass ring.}
		\label{fig:bragg-init-mesh}        
	\end{subfigure}
	\hfill
	\begin{subfigure}[t]{\TWOCOLWIDTH}
		\centering
        \includegraphics[trim = {20cm 8.9cm 20cm 10cm}, clip, width = \ONECOLWIDTH]{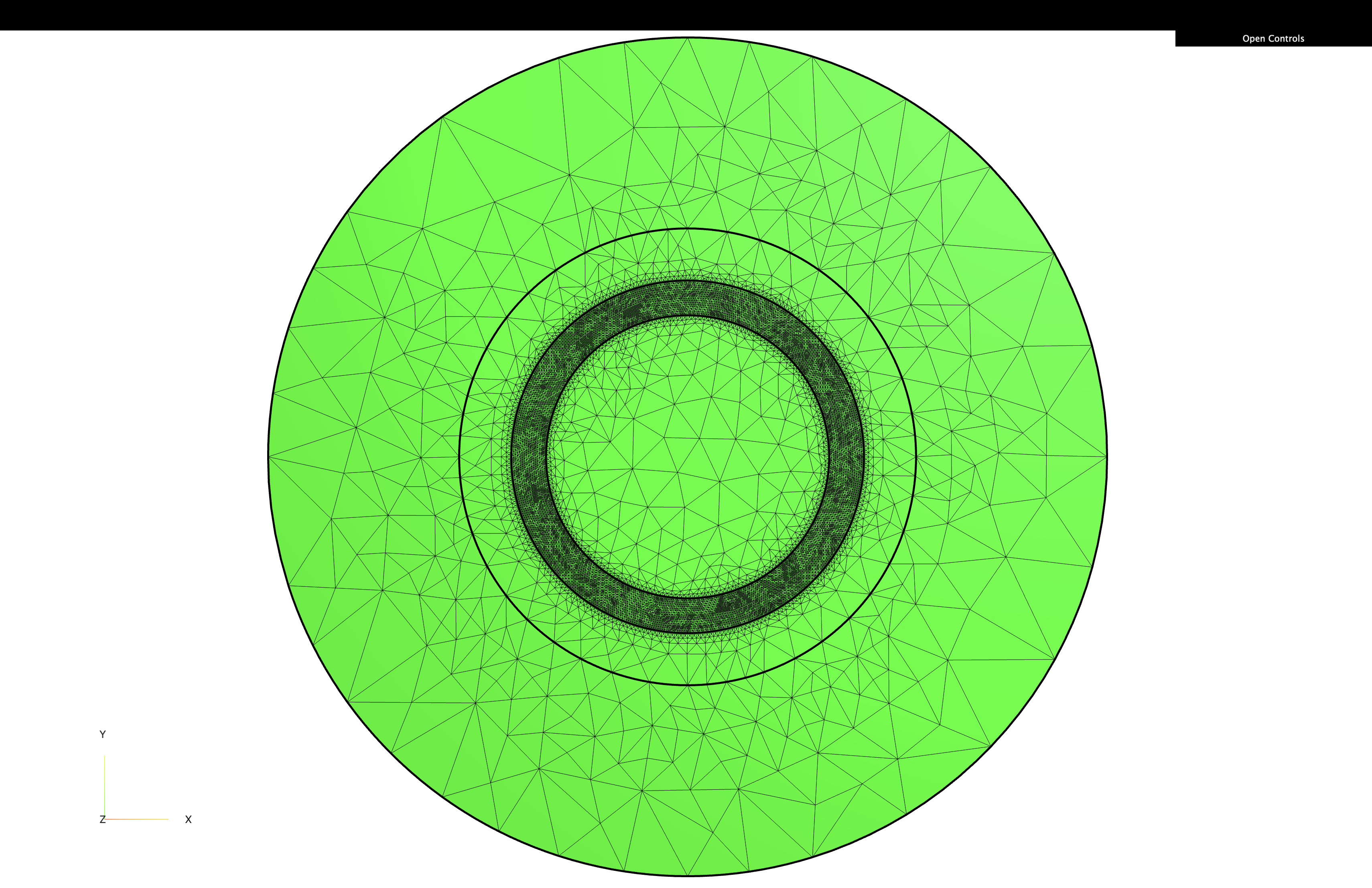}
        \caption{Final mesh for the Bragg fiber after the adaptive algorithm.}
		\label{fig:bragg_mesh25}
	\end{subfigure}
	\vspace{1em}
	\begin{subfigure}[t]{\TWOCOLWIDTH}
		\centering
        \includegraphics[trim = {45cm 22.5cm 45cm 22.5cm}, clip, width = \ONECOLWIDTH]{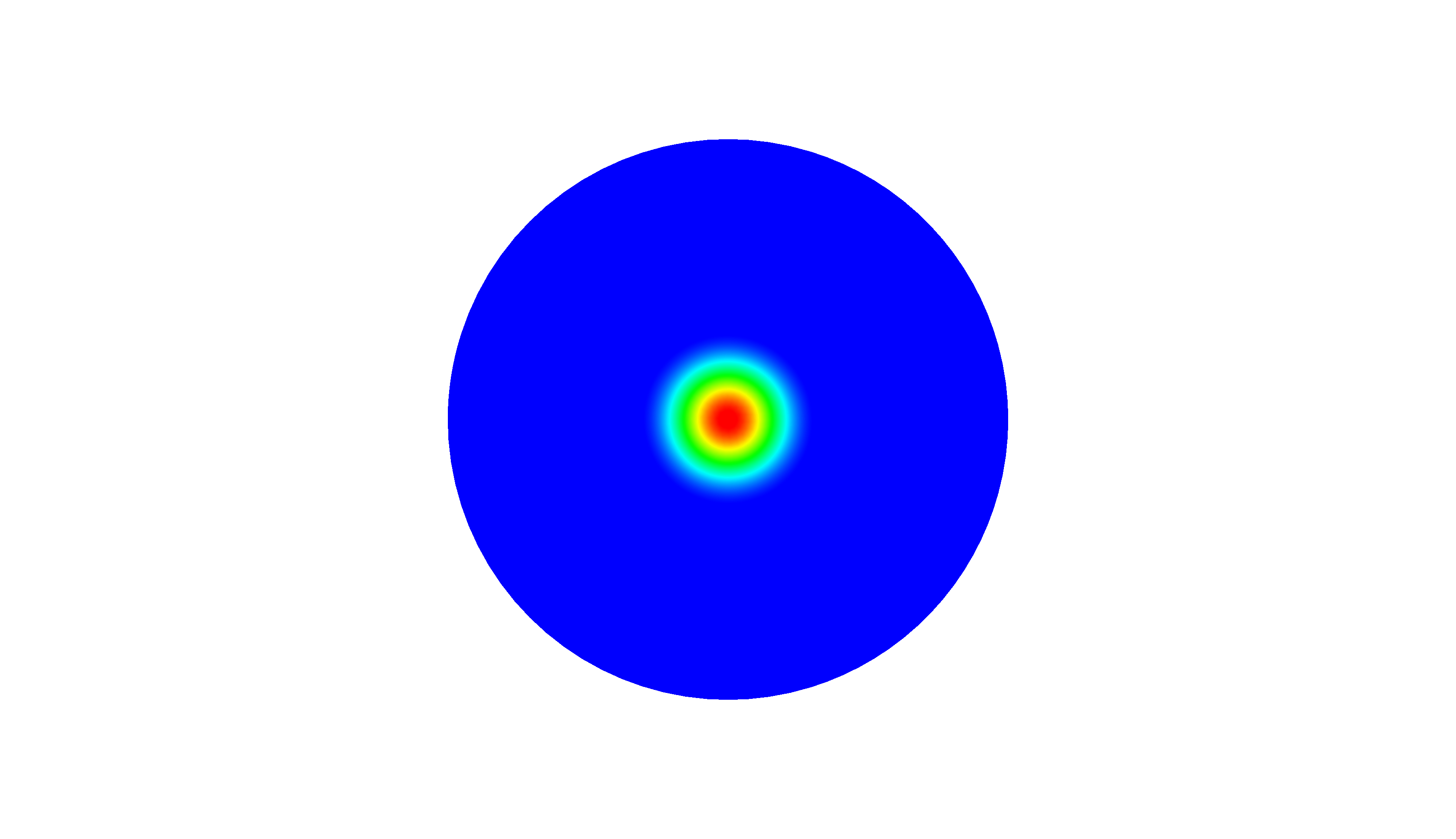}
		\caption{Intensity of the full electric field (right eigenmode) at the final mesh.}
		\label{fig:bragg_intensity}
	\end{subfigure}
	\hfill
	\begin{subfigure}[t]{\TWOCOLWIDTH}    
		\centering
        \includegraphics[trim = {45cm 22.5cm 45cm 22.5cm}, clip, width = \ONECOLWIDTH]{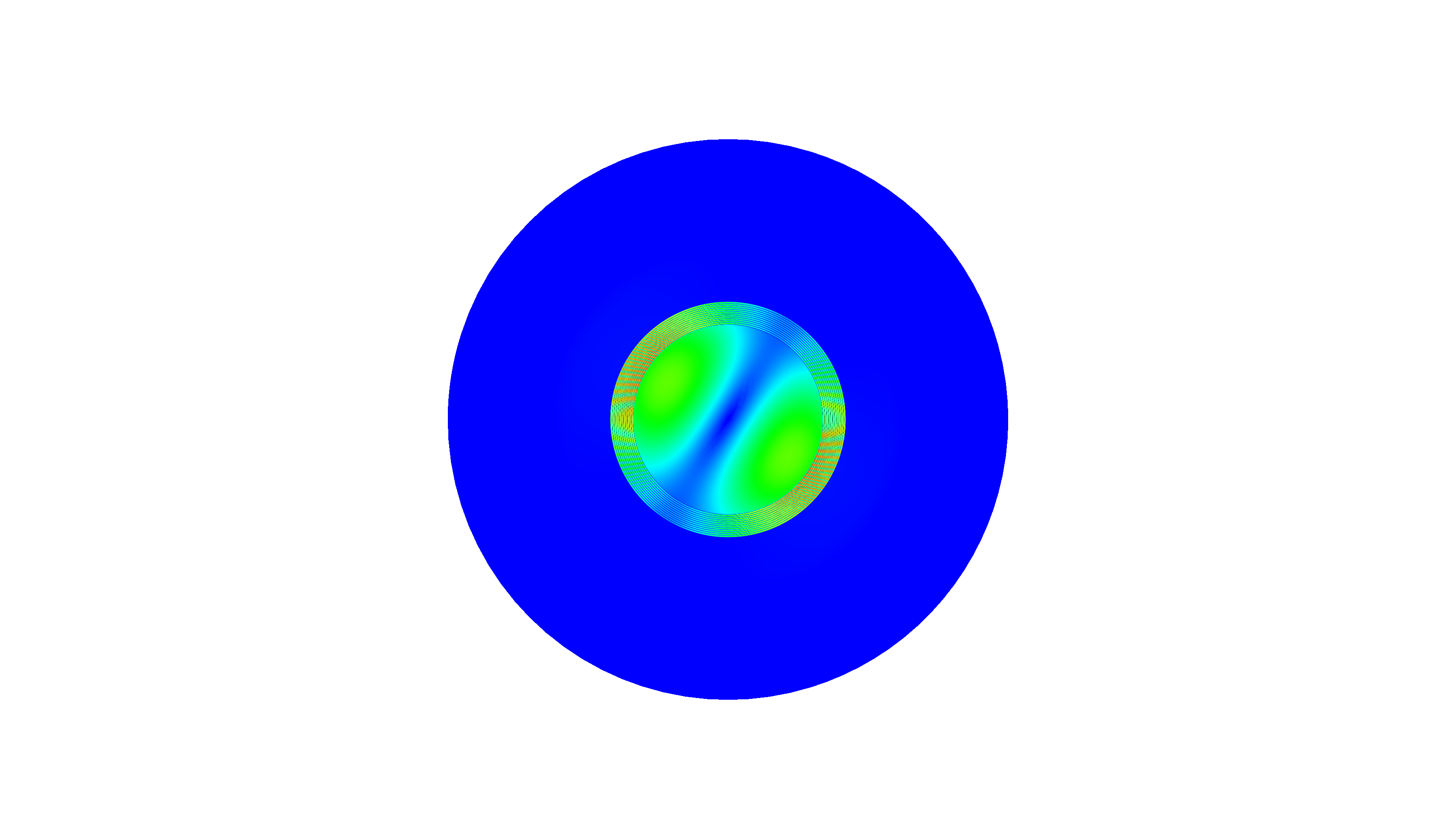}
		\caption{Magnitude of the longitudinal component $(\varphi_h)$ of the eigenmode.}
		\label{fig:bragg_phir}
	\end{subfigure}
	\caption{Results from the adaptive algorithm applied to a Bragg fiber.}\label{fig:bragg-mode}
\end{figure}

In our numerical experiments for this fiber, we use the following settings for parameters (recalling that the parameters were defined in Section~\ref{sec:dwr-maxwell}).  The length scale in nondimensionalization is set to \(L = 1.5 \times 10^{-5}\)~m. 
The operating wavelength is determined by setting \(k = \frac{2 \pi}{1.7} \times 10^{6}~\text{m}^{-1}\). 
The refractive index is a piecewise-constant function whose values are $n_{\rm air} = 1.00027717$ in the air regions and $n_{\rm glass} = 1.43881648$ in the glass ring. 
The following non-dimensionalized radii are used: the air core has radius $r_{\rm core} = 2.7183$; the glass ring starts at $r_{\rm core}$ and has an outer radius $r_{\rm outer} = 3.385$, with a thickness of $t_{\rm ring} = 0.66666667$; the outer air region starts at $r_{\rm outer}$ and extends to $r_{0} = 4.385$, with a thickness of $t_{\rm air} = 1.0$; and the PML starts at $r_{0}$ and extends to $r_{1} = 8.05166666$, with a thickness of $t_{\rm PML} = 3.66666667$. 
The PML strength (recall~\eqref{eq:eta-defn}) is set to $\alpha = 2.0$. 
We set polynomial degree of $p=6$ to define the discretization~\eqref{eq:maxwell-weak-pml-h}. 

The resulting mode at the conclusion of an adaptive iteration is displayed in Figure~\ref{fig:bragg-mode}. 
Clearly, in the initial mesh, visible in Figure~\ref{fig:bragg-init-mesh}, element sizes are more or less the same within the glass ring and in the air regions (with $h \approx  0.1$ in nondimensional units). 
The final mesh produced by the algorithm, given in Figure~\ref{fig:bragg_mesh25}, shows strong local refinement within the glass ring. 
The total electric field intensity at the final iterate (Figure~\ref{fig:bragg_intensity}) reveals a mode highly localized in the air core as expected. In this figure, the fine-scale features are too small in amplitude to be visible. 
However, examining the magnitude of the smaller longitudinal electric field component in Figure~\ref{fig:bragg_phir}, fine-scale ripples in the glass ring are visible. 
These were pointed out in~\cite{2023VanGopGro} from results of semi-analytical computations.  
We see that the algorithm guided the adaptive refinement process to capture these fine-scale features automatically. 
Here and throughout, we use a blue-to-red colormap for mode intensity plots, where blue indicates zero and red indicates the maximal value. 
Color scale is omitted from such plots since the modes are only defined up to a scalar multiple. 

\begin{figure}[htbp]
	\centering
	\begin{tikzpicture}
		\centering
		\begin{loglogaxis}[
				table/col sep=comma,
				xlabel={Number of degrees of freedom at the $\ell$-th iteration},
				ylabel={$d(\varLambda, \varLambda_h)$},
                forget plot style={black!75, dotted},
                ]
            \addplot+ table [x=ndofs, y=hausdorff] {new_tables/efficiency_lvm_adapt.csv};
            %
            %
		\end{loglogaxis}
	\end{tikzpicture}
	\caption{Convergence history of the cluster of discrete eigenvalues to the exact eigenvalue for the Bragg fiber, during the adaptive mesh refinement.}
	\label{fig:bragg-accuracy}
\end{figure}
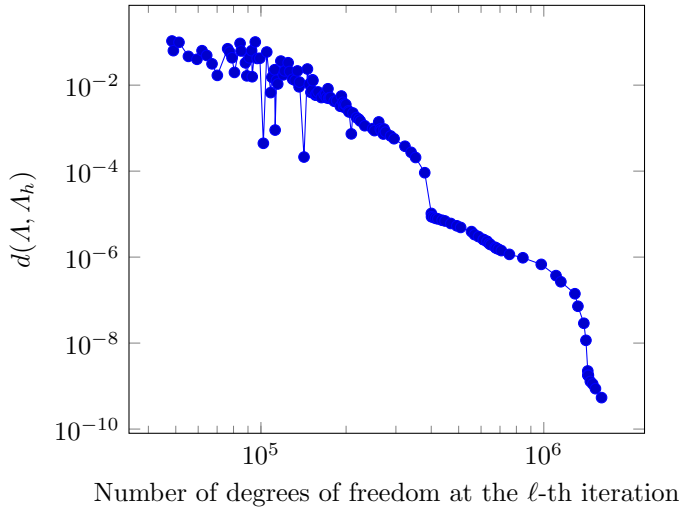

Next, we turn to examine  the accuracy of the computed eigenvalue. The results are in Figure~\ref{fig:bragg-accuracy}. 
The exact nondimensional eigenvalue we aim to approximate, obtained from semi-analytical computations, is  \(Z^2 = 0.80953881+0.00170153i\), and the exact eigenvalue  cluster is the singleton containing it,  $\varLambda = \{ Z^2\}$. 
The discrete eigenvalue cluster $\varLambda_h$ is computed by  providing a circular contour to FEAST centered around the exact eigenvalue. 
When the eigensolver converged, $\varLambda_h$ generally had two elements, matching the multiplicity of the exact eigenvalue. 
Let the Hausdorff distance between two sets $\Upsilon_1$ and $\Upsilon_2$ be denoted by 
\[
	d(\Upsilon_1, \Upsilon_2) = \max( \sup_{\gamma_1 \in \Upsilon_1} \dist(\gamma_1, \Upsilon_2), \sup_{\gamma_2 \in \Upsilon_2} \dist(\gamma_2, \Upsilon_1)) .
\] 
We report the convergence of the discrete eigenvalues in $\varLambda_h$ to $\varLambda$ using the Hausdorff distance in Figure~\ref{fig:bragg-accuracy}. 

\begin{figure}
	\centering
	\begin{tikzpicture}
		\centering
		\begin{loglogaxis}[
			table/col sep=comma,
			xlabel={Number of degrees of freedom at the $\ell$-th iteration},
			ylabel={$|\text{CL}_\ell - \operatorname{CL}|$},
            forget plot style={black!75, dotted},
            ]
            \addplot+ table [x=ndofs, y=error_cl] {new_tables/efficiency_lvm_adapt.csv};

            
		\end{loglogaxis}
	\end{tikzpicture}
	\caption{Convergence history of the confinement loss (see \eqref{eq:CLdef}) for the Bragg fiber, during the adaptive mesh refinement. 
              }
	\label{fig:bragg-CL}
\end{figure}
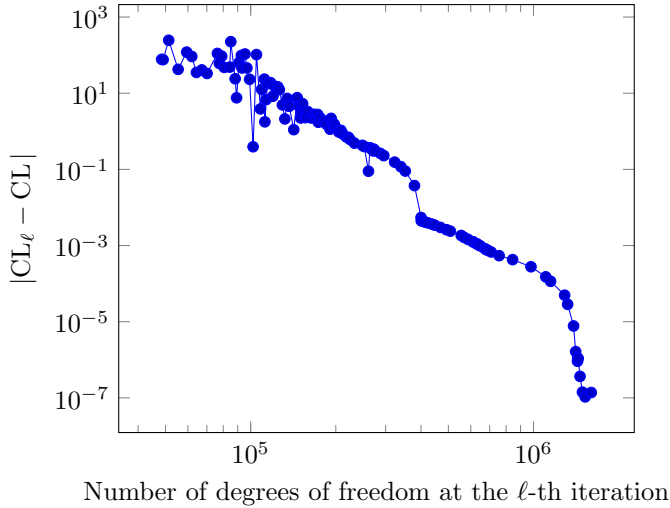

The confinement loss (CL) of a fiber, in decibels per meter (dB/m), is calculated from the propagation constant $\beta$ by
\begin{equation}
	\label{eq:CLdef}
	\text{CL} =
	\text{Im}(\beta) \times \frac{20 }{ \ln(10)} .
\end{equation}
Denote the CL~value obtained using the $\ell$th adaptive eigenvalue iterate by $\text{CL}_\ell$. The difference between this value and the exact confinement loss value calculated from the exact $\beta$ is reported in Figure~\ref{fig:bragg-CL}. 
Excellent agreement is found as adaptive iterations progress. 

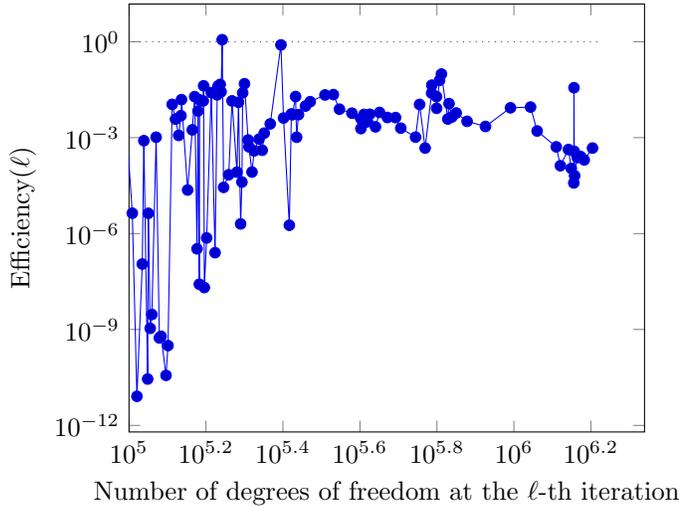
\begin{figure}
	\centering
	\begin{tikzpicture}
		\centering
		\begin{loglogaxis}[
			table/col sep=comma,
			xlabel={Number of degrees of freedom at the $\ell$-th iteration},
            ylabel={Efficiency$(\ell)$}, 
            forget plot style={black!75, dotted},
            xmin=100000, 
            ]
            \addplot [forget plot, black!75, dotted] coordinates {(1, 1) (1650000, 1)};

            \addplot+ table [x=ndofs, y=efficiency_error] {new_tables/efficiency_lvm_adapt.csv};
		\end{loglogaxis}
	\end{tikzpicture}
	\caption{
          \rev{Plot of ratios of eigenvalue error to the estimator (see \eqref{eq:efficiency}) on each mesh during  the adaptive process applied to the Bragg fiber.}
    }
    \label{fig:bragg-efficiency}
\end{figure}

\rev{%
  Finally, we report on how the error estimator compares to the actual eigenvalue error (which can be computed for this example since an exact solution is available).
  Consider the ratio of the the global error  of the computed eigenvalues  to the  global error estimator, i.e., at the $\ell$th
  iteration, let
\begin{equation}
    \label{eq:efficiency}
    \text{Efficiency}({\ell}) =
    \max_{\lambda_h \in \varLambda_h^{(\ell)}} |\lambda_h - \lambda|
    \Bigg/
    \Bigg( \sum_{T \in \oh^{(\ell)}} (\eta_T^{(\ell)})^2 \Bigg)^{1/2},
\end{equation}
where \(\lambda = Z^2\) is the above-mentioned known exact eigenvalue.
These values, for each mesh in the history of adaptive refinements,
are plotted in Figure~\ref{fig:bragg-efficiency}.
As can be seen from
Figure~\ref{fig:bragg-efficiency}, these values oscillates greatly in a
preasymptotic regime, after which they stabilize somewhat,
hovering well below one. On all meshes except one, these values are below
one, indicating that the eigenvalue error is bounded by the computed
error estimator, i.e., the estimator is quite reliable.  However, we also
observe that the efficiency values are not very close to the perfect
value of one.  Such efficiencies are comparable to previous
reports~\cite{2001HeuRan} that lead us to not anticipate  perfect
efficiencies for the DWR estimator,
even for textbook nonselfadjoint eigenproblems simpler than ours.
}%

\section{Three microstructured fibers} \label{sec:three-fibers}

In this section, we report the results obtained using our adaptive algorithm on three different microstructured fiber designs.

\subsection{A hollow-core anti-resonant fiber}\label{sec:arf}

First, we consider the ARF microstructure in Figure~\ref{fig:geo_arf}, motivated by studies in \cite{2013KolKosPryBirPlotDia, 2014Pol} of anti-resonant nodeless tube-lattice fiber. 
This fiber has a hollow-core surrounded by glass capillaries (modeled as) embedded into a glass cladding. 

\begin{figure}[htbp]
	\centering
	\begin{subfigure}[t]{\TWOCOLWIDTH}
		\centering
        \includegraphics[trim = {20cm 8.9cm 20cm 10cm}, clip, width = \ONECOLWIDTH]{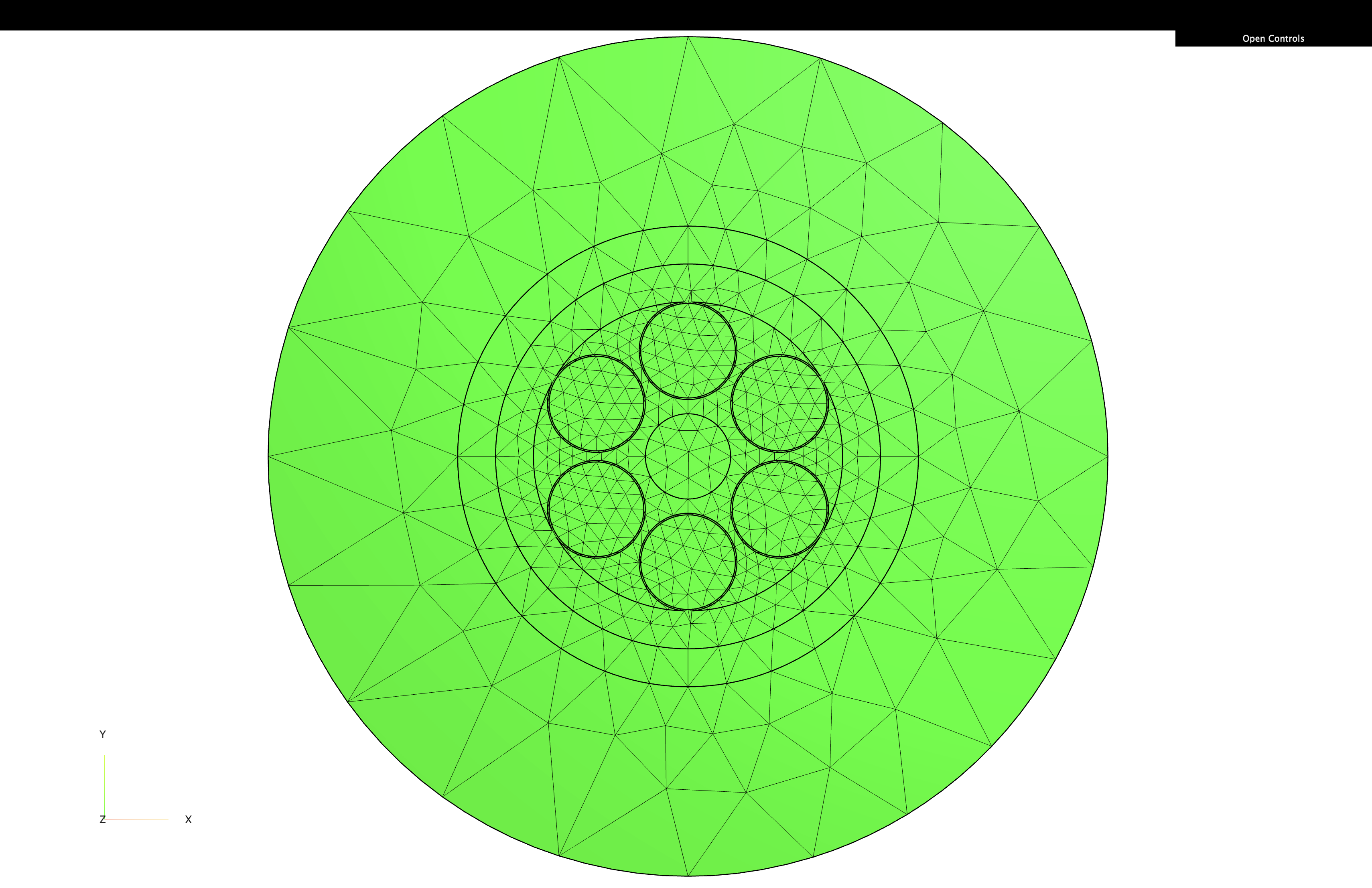}
		\caption{Initial mesh for the ARF geometry.}
		\label{fig:arf_mesh0}
	\end{subfigure}
	\hfill
	\begin{subfigure}[t]{\TWOCOLWIDTH}
		\centering
        \includegraphics[trim = {20cm 8.9cm 20cm 10cm}, clip, width = \ONECOLWIDTH]{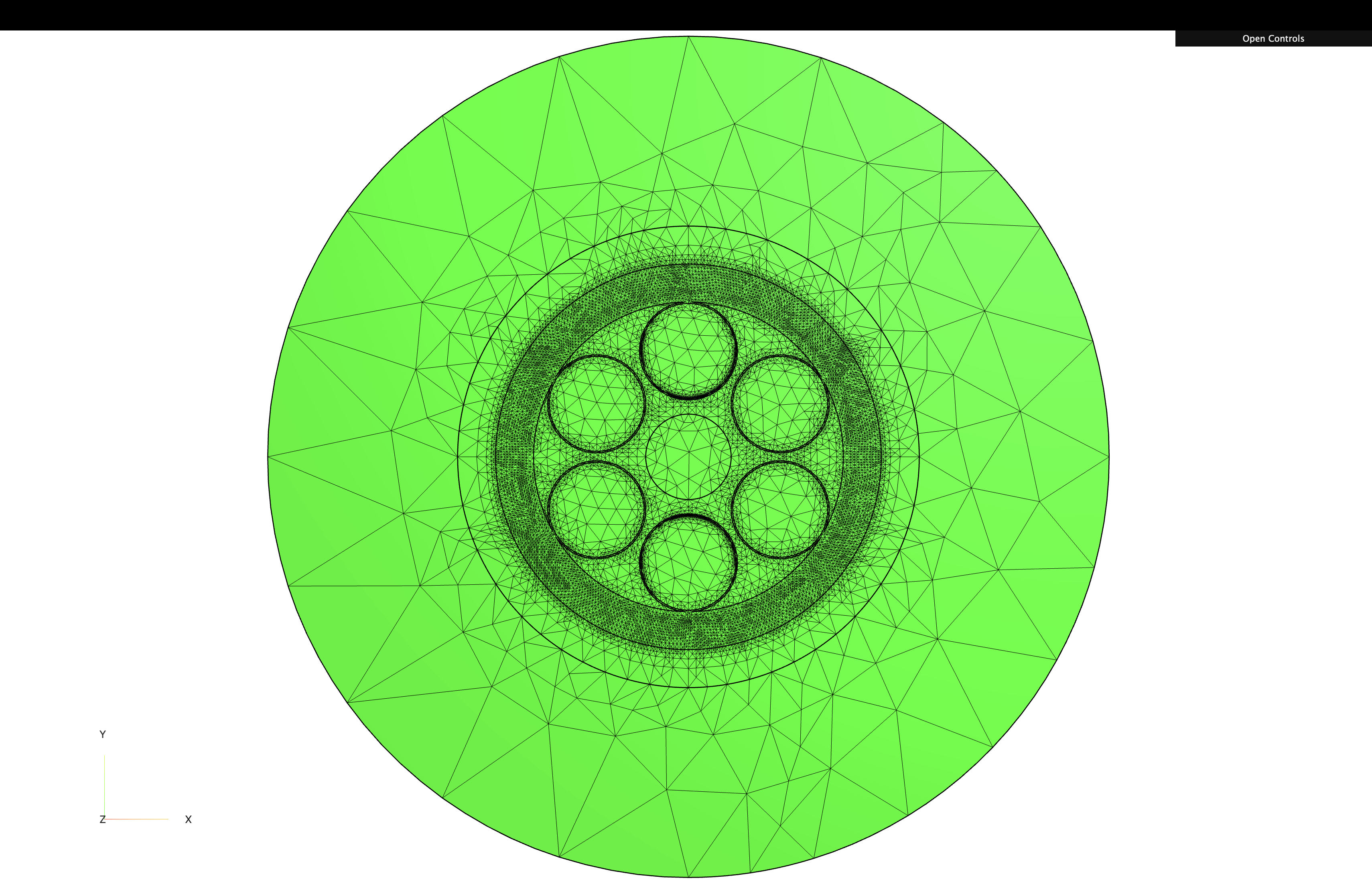}
        \caption{Final mesh for the ARF after the adaptive algorithm}
		\label{fig:arf_mesh}
	\end{subfigure}
	\vspace{1em}
	\begin{subfigure}[t]{\TWOCOLWIDTH}
		\centering
        \includegraphics[trim = {45cm 22.0cm 45cm 22.0cm}, clip, width = \ONECOLWIDTH]{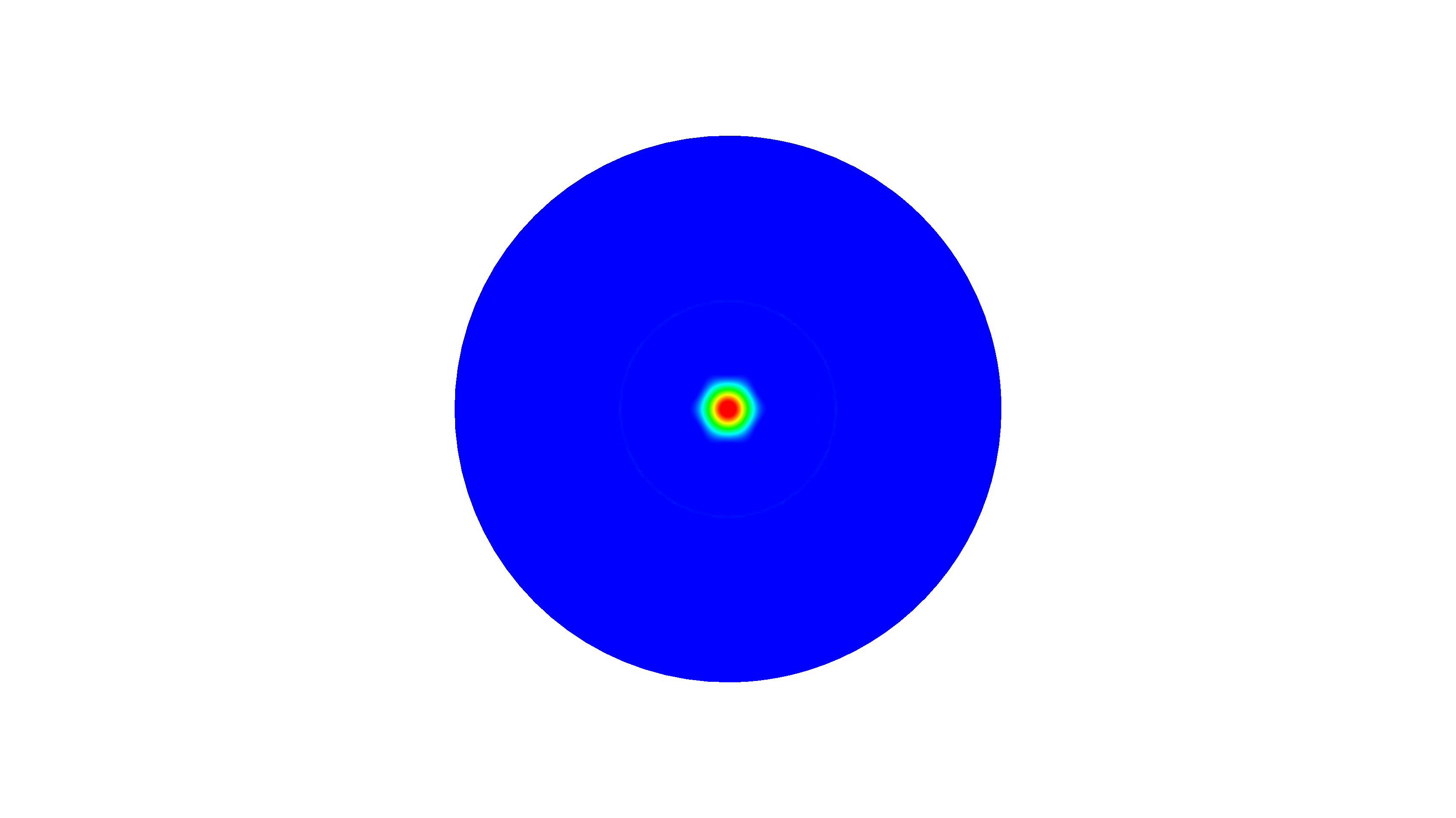}
		\caption{Intensity of the full electric field (right eigenmode) at the final mesh.}
		\label{fig:arf_intensity}
	\end{subfigure}
	\hfill
	\begin{subfigure}[t]{\TWOCOLWIDTH}    
		\centering
        \includegraphics[trim = {45cm 22.5cm 45cm 22.5cm}, clip, width = \ONECOLWIDTH]{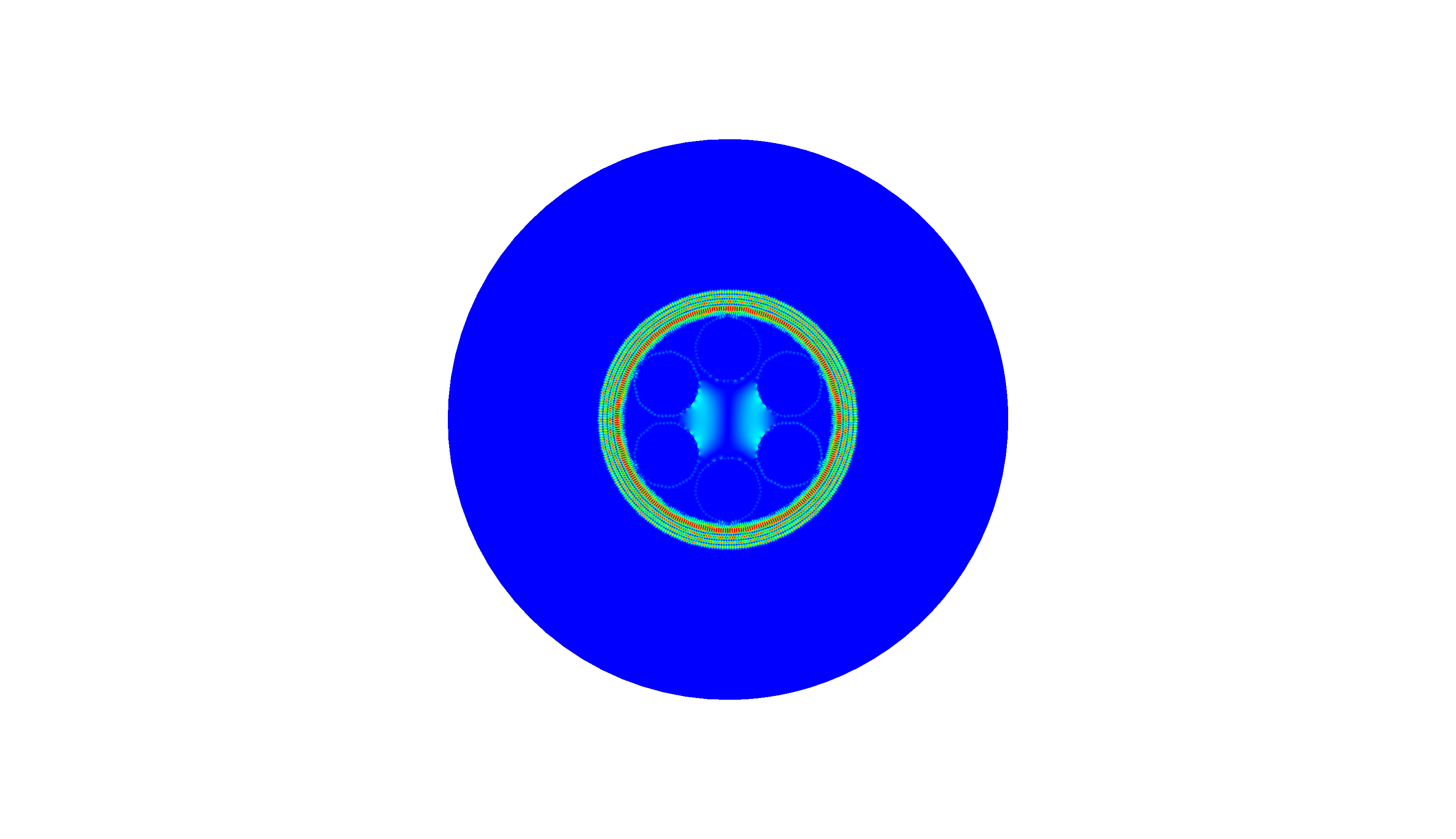}
		\caption{Magnitude of the longitudinal component $(\varphi_h)$ of the eigenmode.}
		\label{fig:arf_phir}
        \end{subfigure}
        \begin{subfigure}[t]{\TWOCOLWIDTH}                
          \begin{tikzpicture}[trim axis left]
            \hspace{-.5cm}
            \begin{loglogaxis}[
              scale=0.8,
              table/col sep=comma,
              ylabel={}, 
              xlabel={Number of d.o.f.s at iteration~$\ell$},
              title={$d(\varLambda_h^{\text{final}},\varLambda_h^{(\ell)})$},
              skip coords between index={0}{1},
              forget plot style={black!75, dotted},
              ]
              \addplot+ table[x=ndofs, y=hausdorff] {new_tables/efficiency_lvm_adapt_arf.csv};


            \end{loglogaxis}
          \end{tikzpicture}
          \caption{History of adaptive convergence. 
          }
          \label{fig:arf-accuracy}
        \end{subfigure}\qquad
        \begin{subfigure}[t]{\TWOCOLWIDTH}
          \begin{tikzpicture}[trim axis right]
            \centering
            \begin{loglogaxis}[
              scale=0.8,
              table/col sep=comma,
              xlabel={Number of d.o.f.s at iteration~$\ell$},
              title={$\text{CL}_\ell$},
              skip coords between index={0}{1},
              forget plot style={black!75, dotted},
              ]
              \addplot+ table [x=ndofs, y=cl] {new_tables/efficiency_lvm_adapt_arf.csv};
            \end{loglogaxis}
          \end{tikzpicture}
          \caption{History of computed loss values. 
          }
          \label{fig:arf-cl}
        \end{subfigure}
	\caption{Results from the adaptive algorithm applied to ARF.}
	\label{fig:arf-results}
\end{figure}

The parameters used  for the simulation are as follows. 
Length scale is \(L = 1.5  \times  10^{-5}\)~m; operating wavelength is determined by \(k = \frac{2 \pi}{1.8}  \times  10^{6}~\text{m}^{-1}\); the refractive index is $n_{\rm air} = 1.00027717$ in the air core, and $n_{\rm glass} = 1.44087350$ in the glass ring. 
The following non-dimensionalized radii and thicknesses are used: the core of the fiber is a circle of radius $r_{\rm core} = 1.0$; the capillary is defined between an inner radius $r_{\rm cap,\ i} = 2.71833333$ and an outer radius $r_{\rm cap,\ o} = 3.385$, with a thickness of $t_{\rm cap} = 0.66666667$; the cladding has an inner radius $r_{\rm clad,\ i} = 3.385$ and an outer radius $r_{\rm clad,\ o} = 4.05166667$, the thickness of the glass jacket (cladding) is $t_{\rm clad} = 0.66666667$; the embedding of the capillary into the cladding is given by \(e = 0.05952380\), which represents the fraction of the capillary thickness that is submerged in the cladding; the distance between the capillaries is \(d = 0.13999999\); the outer layer of air starts at $r_{\rm clad,\ o}$ and ends at $r_{0} = 4.05166667$, and has a thickness of $t_{\rm air} = 0.66666667$; the PML region extends from $r_{0}$ to $r_{1} = 7.385$, with a thickness of $t_{\rm PML} = 3.33333333$. 
The PML strength (implemented through relation~\eqref{eq:eta-defn}) is $\alpha = 2.0$. 
The algorithm was run with a search region to capture the expected core-localized fundamental mode of this fiber. 

The results are portrayed in the images of Figure~\ref{fig:arf-results}. 
The initial mesh (Figure~\ref{fig:arf_mesh0}) is relatively coarse, with small elements used only to conform to the geometry of the thin glass capillaries and their melding into the cladding. 
We make the elements in the core region smaller (with \(0.2 \lesssim h \lesssim 0.5\)), to capture the expected core-localized fundamental mode of this fiber and verify the saturation assumption. 
The final mesh produced by the adaptivity iteration (Figure~\ref{fig:arf_mesh}) is characterized by unexpectedly strong local refinement in the outer glass cladding where there are no tiny geometrical features to be resolved. 
This is explained by Figure~\ref{fig:arf_phir} where we see fine-scale ripples in the high-index cladding for the longitudinal electric field component. 
The total intensity (Figure~\ref{fig:arf_intensity}) is characterized by a localized profile in the air core. 
However, because the mode intensity in the core region is orders of magnitude larger than its fine-scale features in the cladding region, one may not realize such cladding ripples exist and are important for accurately resolving the mode loss. 
The cladding oscillations of the mode intensity are visible when plotting the smaller longitudinal vector component of the mode. 
The importance of capturing these fine scale features was first pointed out in~\cite{2023VanGopGro}, where an expert ``informed mesh'' was manually created with enough elements to capture the fine-scale oscillations. 
Here, we see that the automatic adaptive process also points to the necessity of refinements in the same region.

The apparent convergence history during the adaptive process is displayed in Figure~\ref{fig:arf-accuracy}. 
Since we do not have exact eigenvalues for this case, we cannot calculate the eigenvalue errors in each iteration. 
Instead, we compare the Hausdorff distance between the final output $\varLambda_h^{\text{final}}$ (at the last iteration) with the discrete eigenvalue cluster $\varLambda_h^{(\ell)}$ (with two elements converging to the same number) found at the $\ell$th adaptive iteration. 
The curve exhibits an overall decreasing trend, even if not a monotonic decrease. 
The CL values at each iteration are compared to the CL value calculated from the final eigenvalue iterate in Figure~\ref{fig:arf-cl}, which illustrates convergence of the CL approximations.

\subsection{A nested anti-resonant nodeless fiber}\label{sec:nanf}

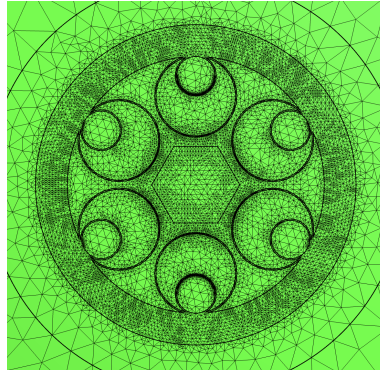
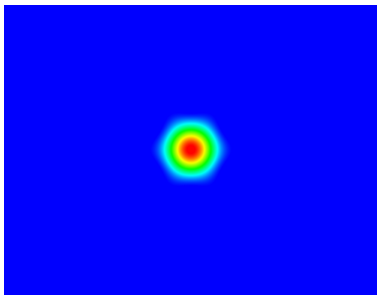
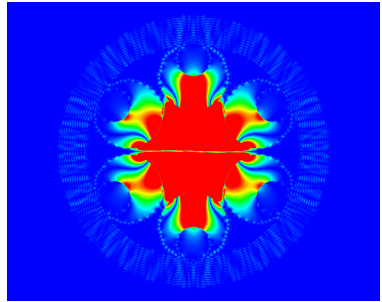
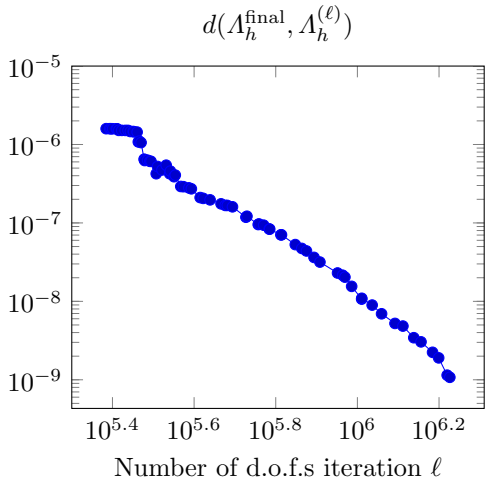
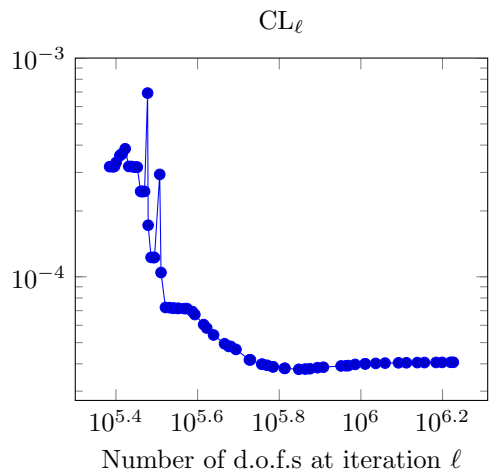
\begin{figure}
	\centering
	\begin{subfigure}[t]{\TWOCOLWIDTH}
		\centering
        \includegraphics[trim = {20cm 8.9cm 20cm 10cm}, clip, width = \ONECOLWIDTH]{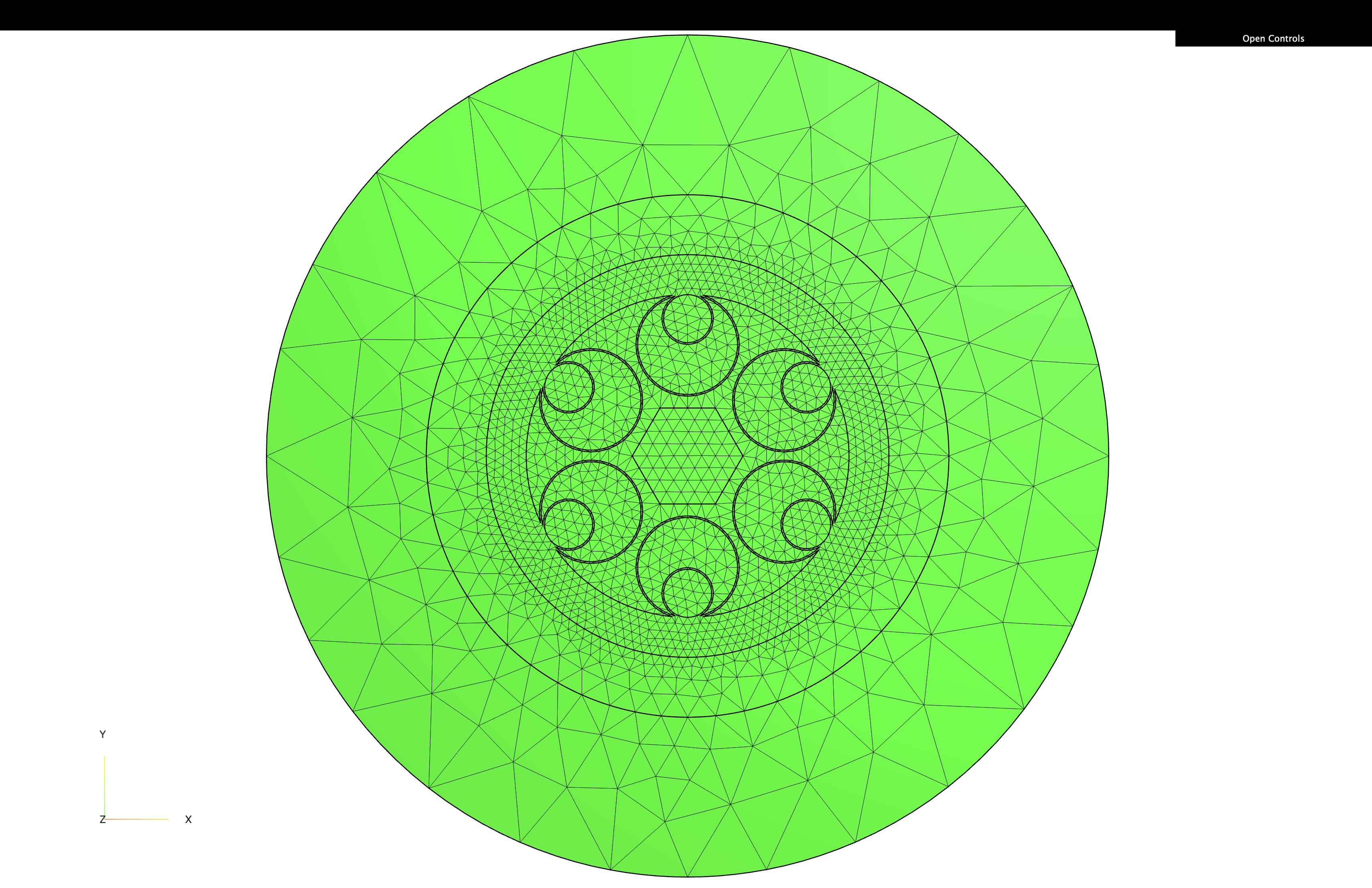}
        \caption{Initial mesh for the NANF fiber.}
		\label{fig:nanf_mesh0}
	\end{subfigure}
	\hfill
	\begin{subfigure}[t]{\TWOCOLWIDTH}
		\centering
        \includegraphics[trim = {20cm 8.9cm 20cm 10cm}, clip, width = \ONECOLWIDTH]{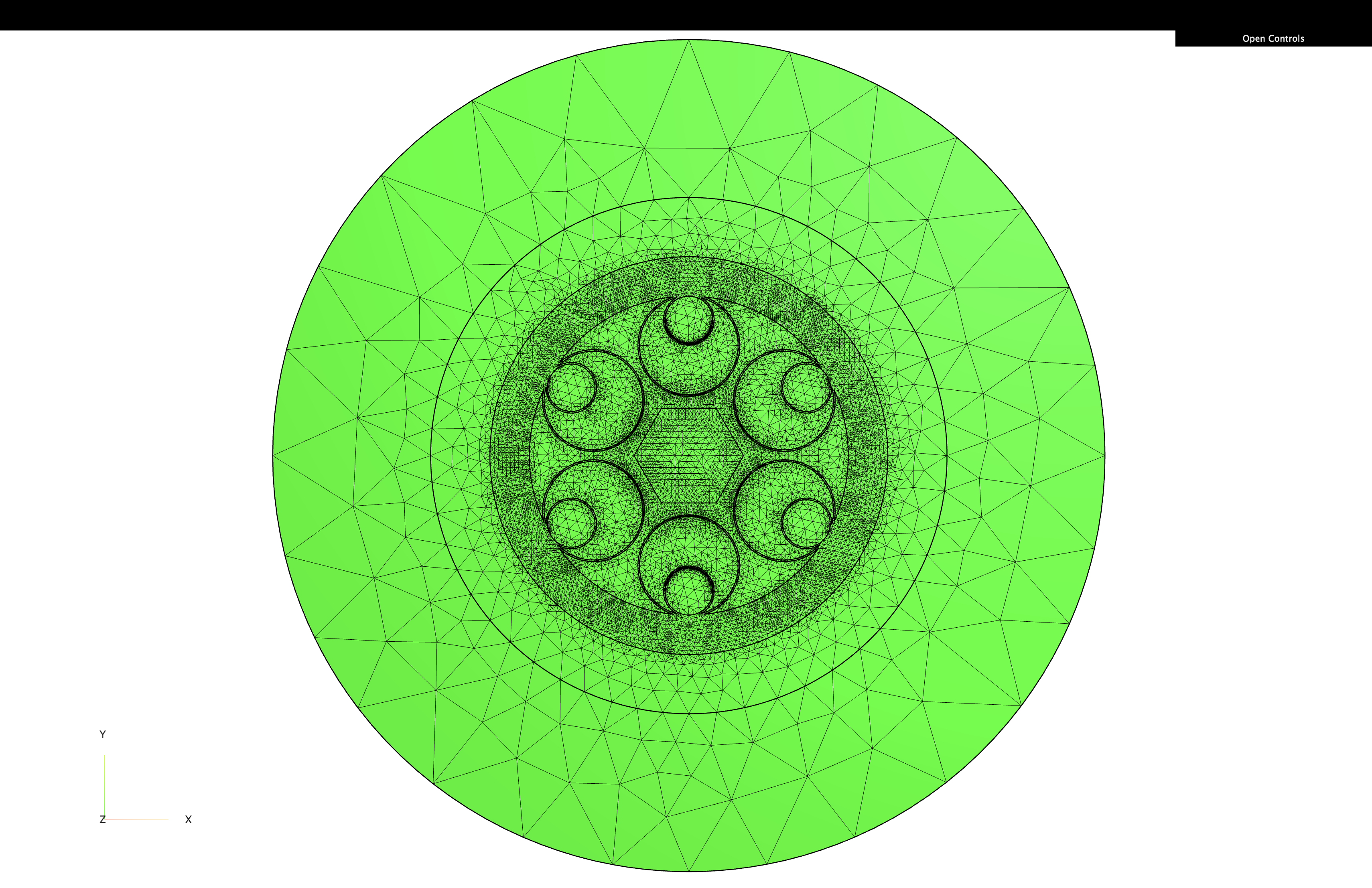}
        \caption{Final mesh for the NANF fiber after the adaptive algorithm.}
		\label{fig:nanf_mesh34}
	\end{subfigure}
	\vspace{1em}
	\begin{subfigure}[b]{\TWOCOLWIDTH}
		\centering
        \includegraphics[trim = {45cm 22.0cm 45cm 22.0cm}, clip, width = \ONECOLWIDTH]{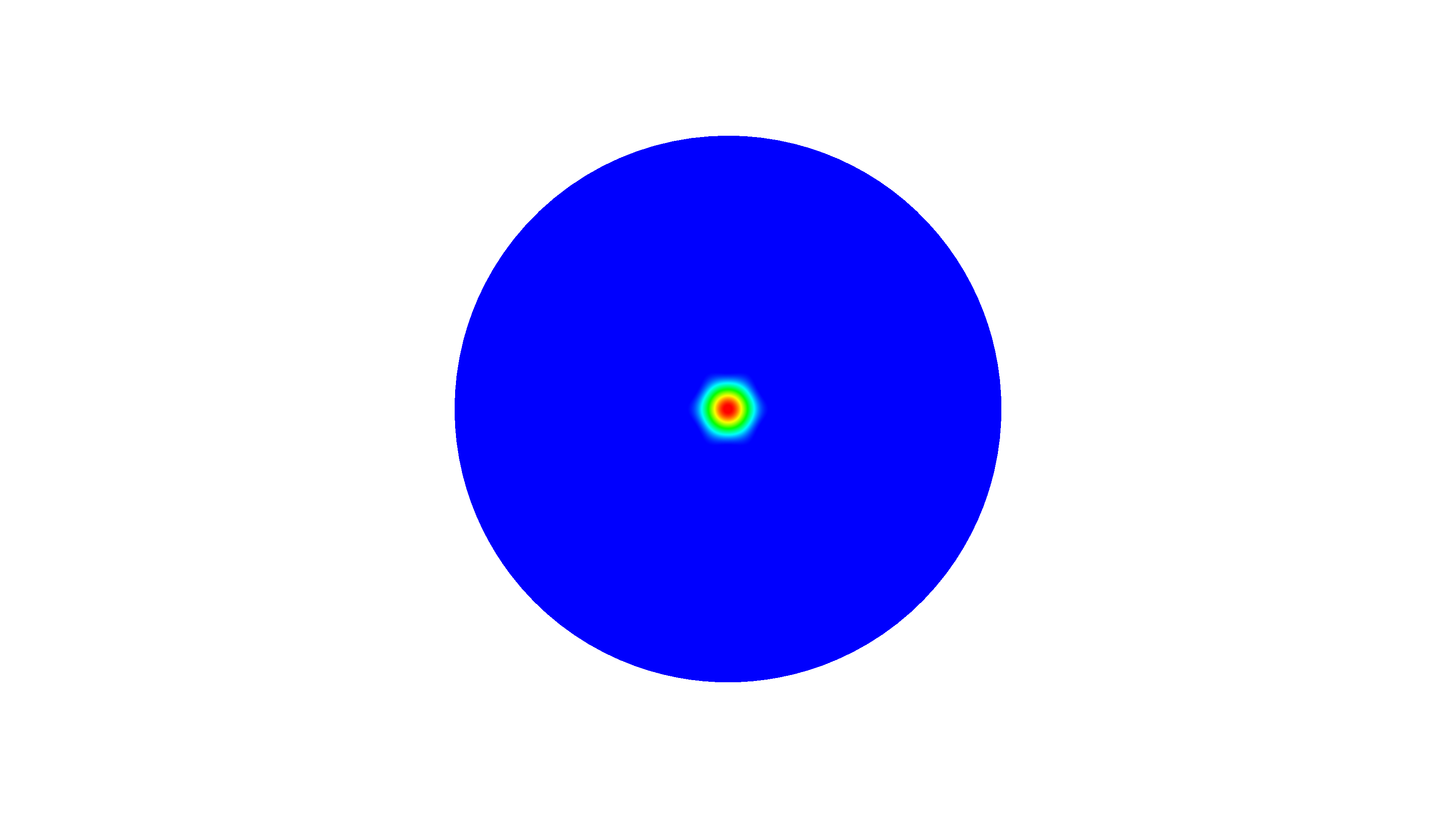}
		\caption{Intensity of the full electric field (right eigenmode) at the final mesh.}
		\label{fig:nanf_intensity}
	\end{subfigure}
	\hfill
	\begin{subfigure}[b]{\TWOCOLWIDTH}    
		\centering
        \includegraphics[trim = {45cm 22.5cm 45cm 22.5cm}, clip, width = \ONECOLWIDTH]{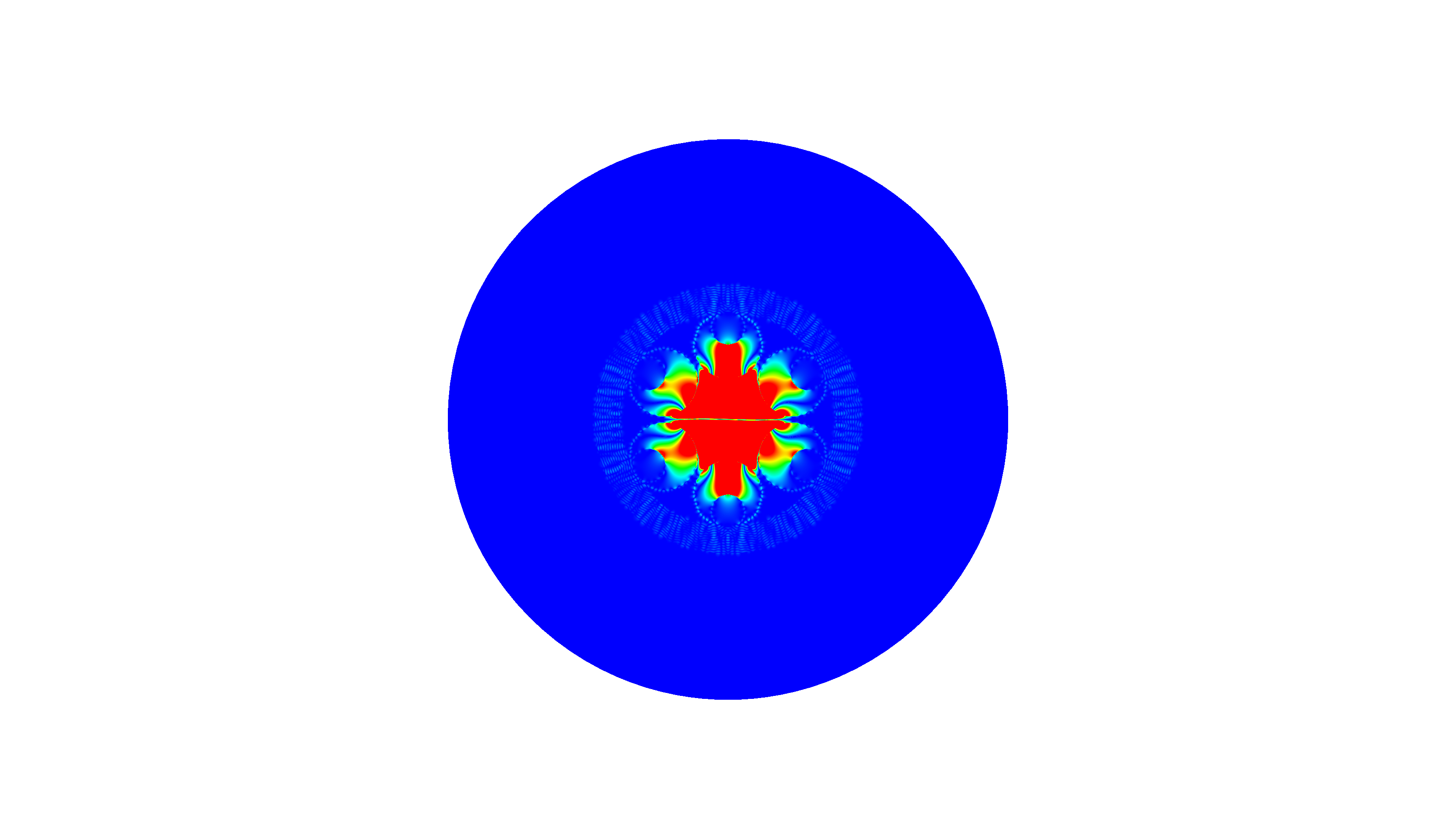}
		\caption{Magnitude of the longitudinal component $(\varphi_h)$ of the eigenmode.}
		\label{fig:nanf_phir}
       \end{subfigure}
       \begin{subfigure}[t]{\TWOCOLWIDTH}
         \begin{tikzpicture}[trim axis left]
           \hspace{-.5cm},
           \begin{loglogaxis} [
             scale=0.8,
             table/col sep=comma,
             xlabel={Number of d.o.f.s  iteration  $\ell$}, 
             title={$d(\varLambda_h^{\text{final}},\varLambda_h^{(\ell)})$},
             skip coords between index={0}{1},
             ymax=1.0e-5,
             forget plot style={black!75, dotted},
             ]
             \addplot+ table[x=ndofs, y=hausdorff] {new_tables/efficiency_lvm_adapt_nanf.csv};


           \end{loglogaxis}
         \end{tikzpicture}
         \caption{History of adaptive convergence.} 
         \label{fig:nanf-accuracy}

       \end{subfigure}
       \quad 
       \begin{subfigure}[t]{\TWOCOLWIDTH}
         \begin{tikzpicture}[trim axis right]
           \centering
           \begin{loglogaxis}[
             scale=0.8,
             table/col sep=comma,
             xlabel={Number of d.o.f.s at iteration $\ell$},
             title={$\text{CL}_\ell$},
             ymax=1.0e-3,
             skip coords between index={0}{1},
             forget plot style={black!75, dotted},
             ]
             \addplot+ table [
             x=ndofs,
             y=cl,
             ] {new_tables/efficiency_lvm_adapt_nanf.csv};

           \end{loglogaxis}
         \end{tikzpicture}
         \caption{History of  computed loss values.} 
         \label{fig:nanf-cl}
       \end{subfigure}
	\caption{Results from the adaptive algorithm applied to NANF.}
	\label{fig:nanf-all}
\end{figure}

Next, we consider the NANF microstructure in Figure~\ref{fig:geo_nanf}, motivated by the studies in \cite{2014Pol, 2018BraHayCheJasSanSlaFokBawSakDav}. 
Its main difference with the ARF is that the cladding has a nested structure of capillaries. 
Since this structure is close to the previously considered ARF, we anticipate that fine-scale mode features of a similar nature  to arise here. 

We employ the following parameters for the simulation. 
The length scale is \(L = 1.5  \times  10^{-5}\)~m; the operating wavelength is determined by \(k = \frac{2 \pi}{1.8} \times 10^{6}~\text{m}^{-1}\); the refractive index is $n_{\rm air} = 1.00027717$ in the air regions, and $n_{\rm glass} = 1.44087350$ in the glass ring. 
We utilize the following non-dimensionalized radii and thicknesses: the core of the fiber is a circle of radius $r_{\rm core} = 1.0$; the outer capillary is determined by an outer radius $r_{\rm cap,\ o} = 0.832$, and a thickness of $t_{\rm cap,\ o} = 0.028$; the inner capillary is defined by an inner radius $r_{\rm cap,\ i} = 0.4$, and a thickness of $t_{\rm cap,\ i} = 0.028$; the cladding has an inner radius $r_{\rm clad} = 1.33333333$, an outer radius $r_{\rm buffer} = 2.0$, and a thickness $t_{\rm clad} = 0.66666667$; the glass ring has an inner radius $r_{\rm inner} = 2.692$, an outer radius $r_{\rm outer} = 4.35866667$, and a thickness $t_{\rm clad} = 1.66666667$; and the PML region starts at $r_{\rm outer}$ and ends on $r_{\rm PML} = 7.02533333$. 
We consider a PML strength (see~\eqref{eq:eta-defn}) of value $\alpha = 2.0$.

The results  in Figure~\ref{fig:nanf-all} confirm that fine-scale modal features exist for this fiber as well.
The initial mesh (Figure~\ref{fig:nanf_mesh0}) is relatively coarse, with small elements used only to conform to the geometry of the thin glass capillaries and their melding into the cladding (as in the ARF case). 
We make the elements in the core region smaller (with \(0.2 \lesssim h \lesssim 0.5\)), to capture the expected core-localized fundamental mode of this fiber and verify the saturation assumption. 
Again, the adaptive meshing algorithm found (in Figure~\ref{fig:nanf_mesh34}) that significant refinements were needed in the outer glass cladding; however, it also determined that some refinements were warranted in the hollow-core near larger capillaries that bound the core. 
The captured fine-scale modal features are visible in Figure~\ref{fig:nanf_phir}. 
The convergence of the eigenvalues and CL values  during the adaptive process are illustrated in Figures~\ref{fig:nanf-accuracy} and~\ref{fig:nanf-cl}.

\subsection{A photonic bandgap fiber}\label{sec:pbg}

Finally, we consider the PBG microstructure design in Figure~\ref{fig:geo_pbg},  previously studied in~\cite{2013PolPetRic, 1999CreManKniBirRusRobAll, 2003LitDunUsnEggWhiMcPdeSte}. 
The lattice arrangement of the dielectric rods in the fiber can be modified to change the guiding properties of the fiber. 
We base our model parameters  on the descriptions in~\cite{2003LitDunUsnEggWhiMcPdeSte}. 

\begin{figure}[htbp]
	\centering
	\begin{subfigure}[t]{\TWOCOLWIDTH}
		\centering
        \includegraphics[trim = {21.5cm 10cm 21.5cm 11cm}, clip, width = \ONECOLWIDTH]{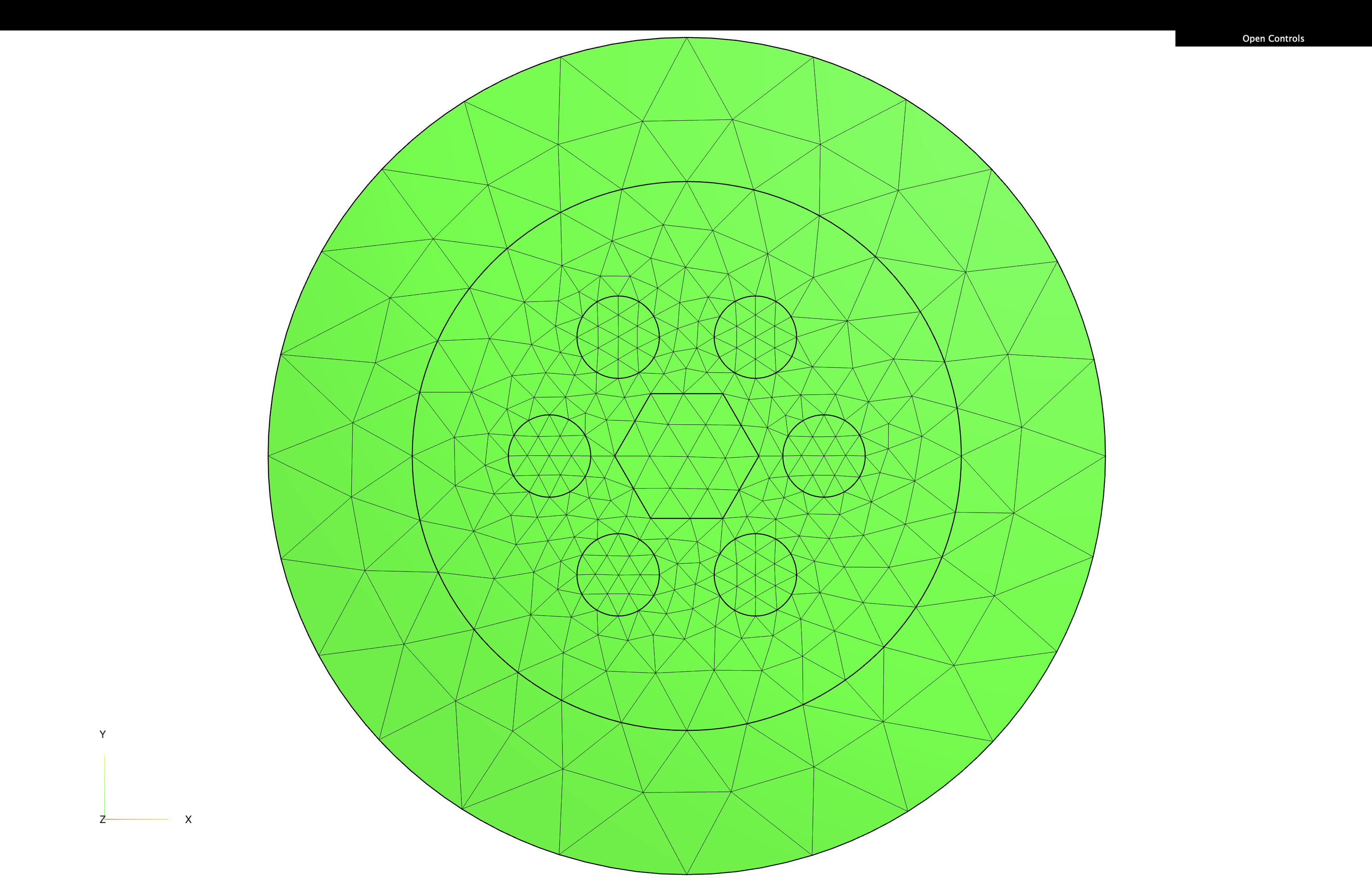}
		\caption{Initial mesh for the PBG fiber, zoomed in near the core}
		\label{fig:pbg_mesh0}
	\end{subfigure}
	\hfill
	\begin{subfigure}[t]{\TWOCOLWIDTH}
		\centering
        \includegraphics[trim = {21.5cm 10cm 21.5cm 11cm}, clip, width = \ONECOLWIDTH]{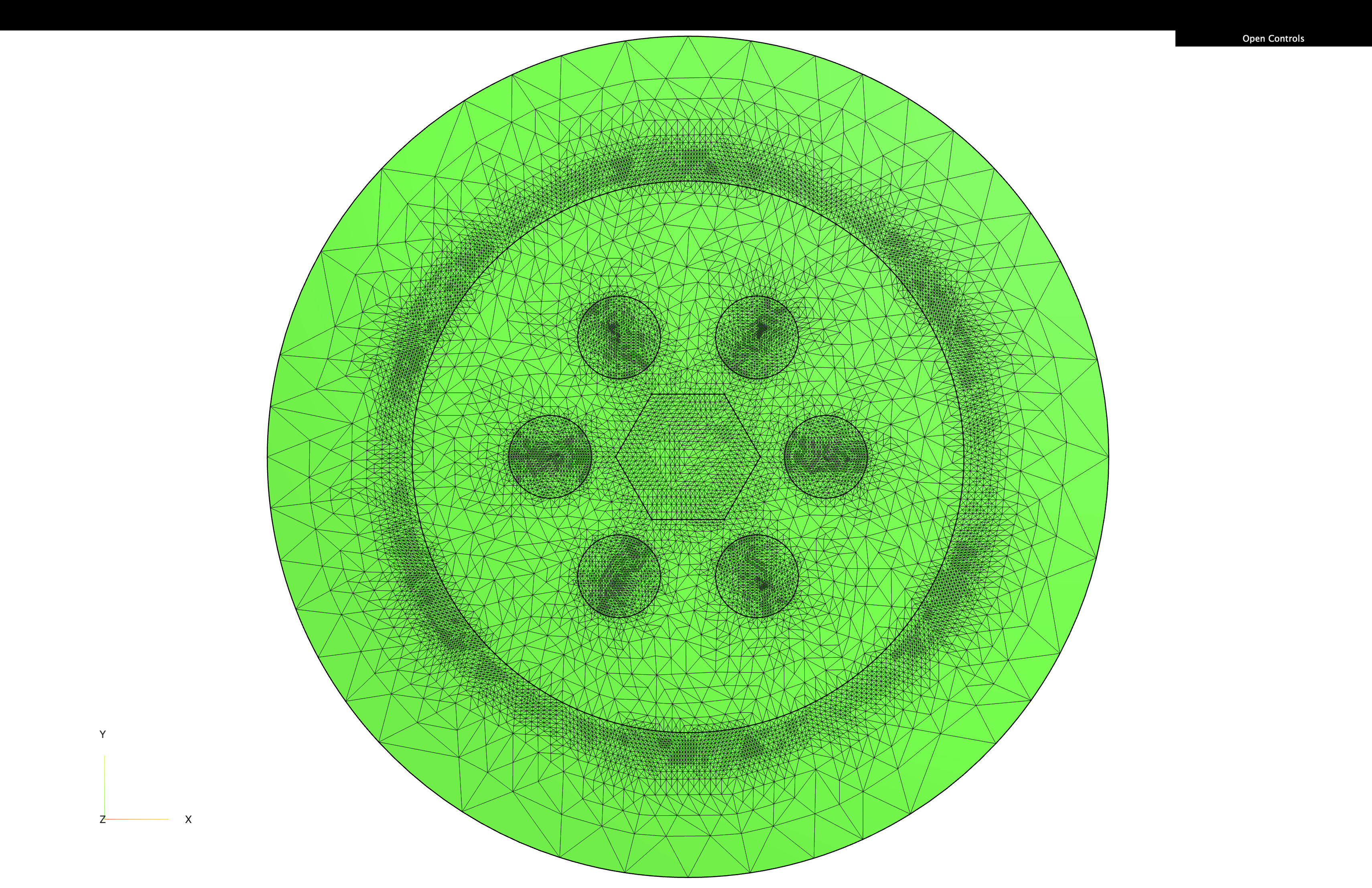}
        \caption{Final mesh for the PBG fiber after the adaptive algorithm.}
		\label{fig:pbg_mesh13}
	\end{subfigure}
	\vspace{1em}
	\begin{subfigure}[t]{\TWOCOLWIDTH}
		\centering
        \includegraphics[trim = {45cm 22.0cm 45cm 22.0cm}, clip, width = \ONECOLWIDTH]{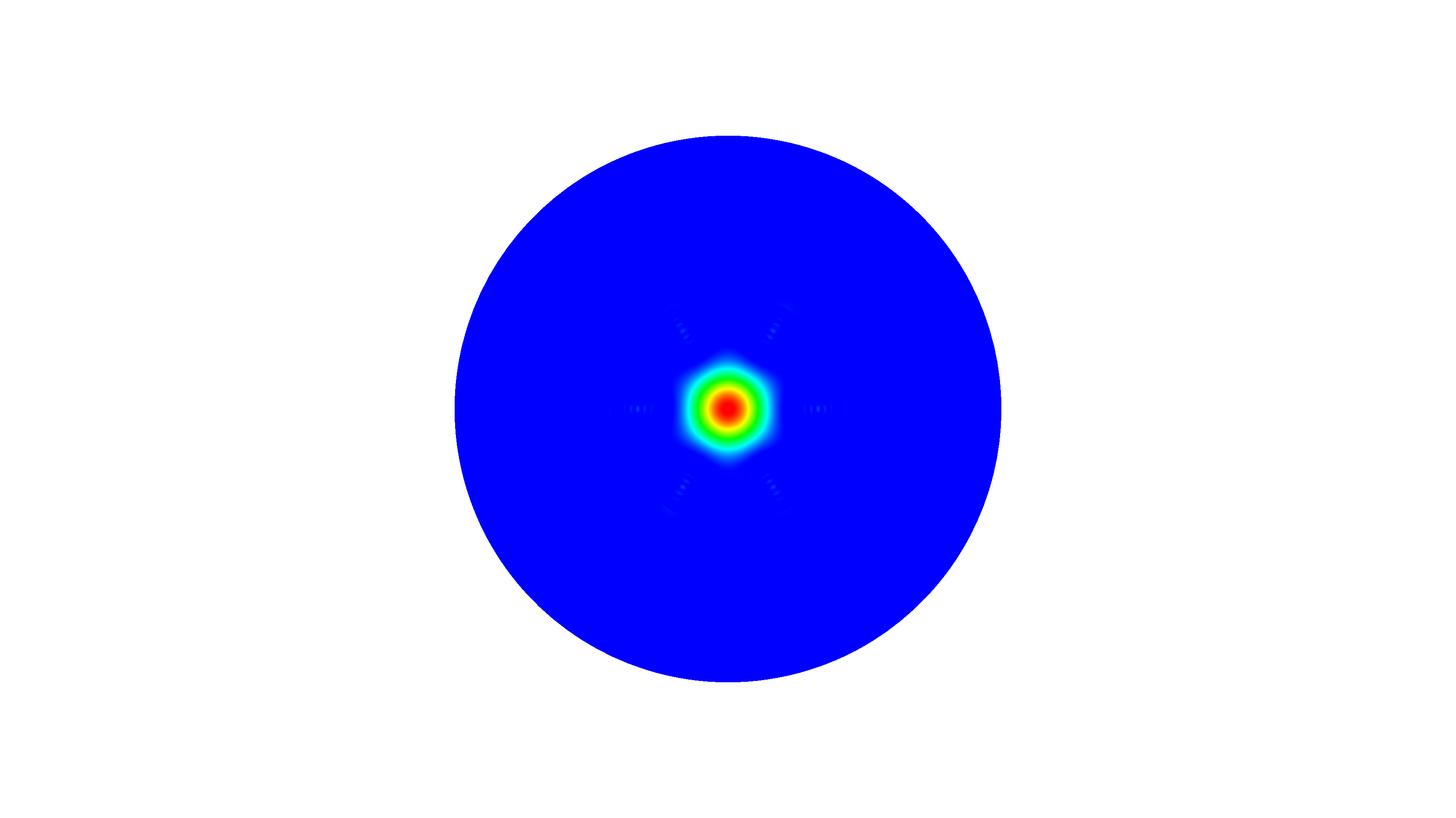}
		\caption{Intensity of the full electric field (right eigenmode) at the final mesh.}
		\label{fig:pbg_intensity}
	\end{subfigure}
	\hfill
	\begin{subfigure}[t]{\TWOCOLWIDTH}    
		\centering
        \includegraphics[trim = {45cm 22.5cm 45cm 22.5cm}, clip, width = \ONECOLWIDTH]{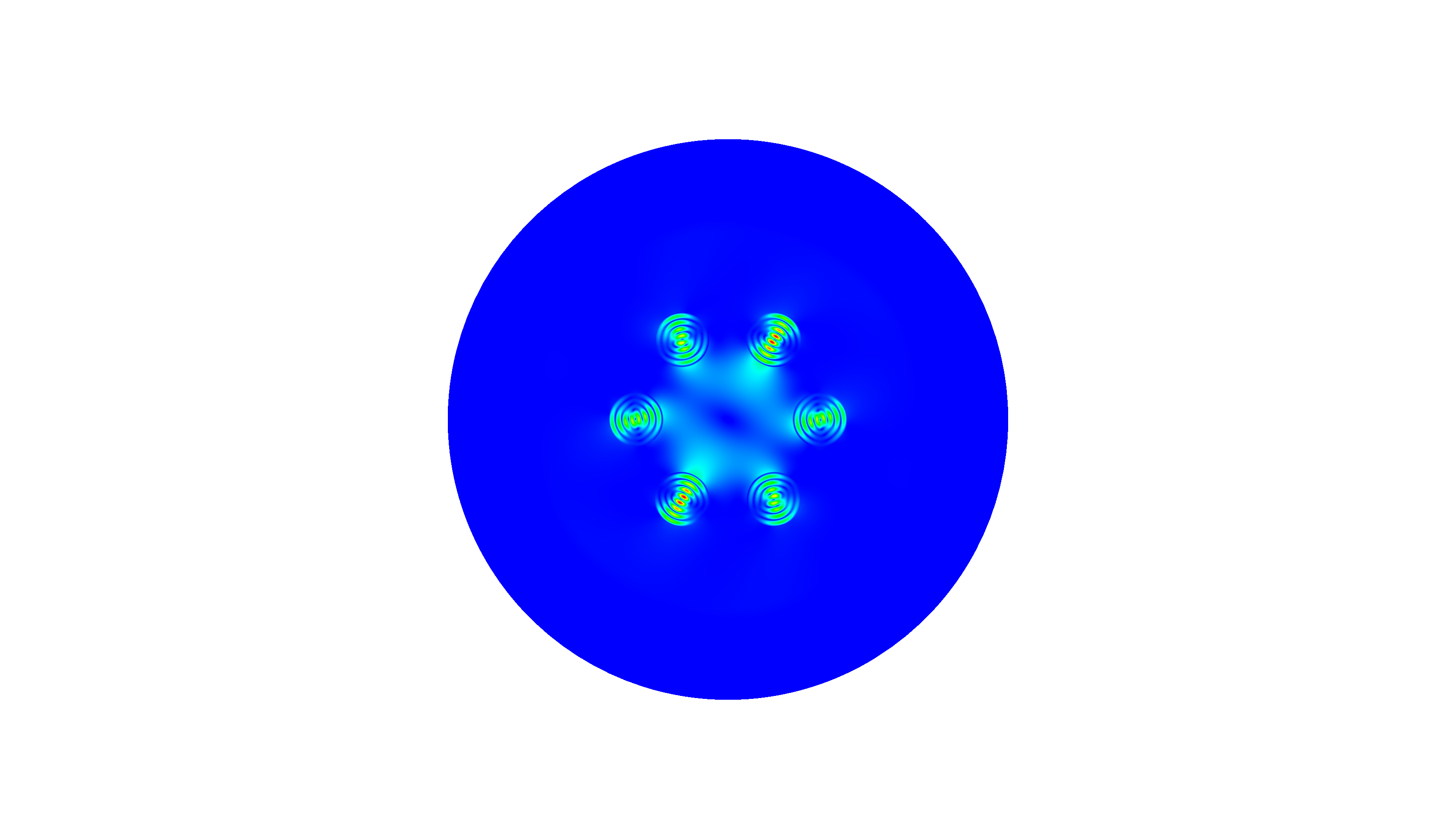}
		\caption{Magnitude of the longitudinal component $(\varphi_h)$ of the eigenmode.}
		\label{fig:pbg_phi}
        \end{subfigure}
        \begin{subfigure}[t]{\TWOCOLWIDTH}    
          \centering
          \begin{tikzpicture}[trim axis left]
            \hspace{-.5cm}
            \begin{loglogaxis}[
              scale=0.8,
              table/col sep=comma,
              xlabel={Number of d.o.f.s at iteration~$\ell$},
              title={$d(\varLambda_h^{\text{final}},\varLambda_h^{(\ell)})$},
              forget plot style={black!75, dotted},
              ]
              \addplot+ table[x=ndofs, y=hausdorff]
              {new_tables/efficiency_lvm_adapt_pbg.csv};
              
            \end{loglogaxis}
	\end{tikzpicture}
	\caption{History of adaptive convergence.}
	\label{fig:pbg-accuracy}
        \end{subfigure}
        \quad 
        \begin{subfigure}[t]{\TWOCOLWIDTH}    
          \centering
          \begin{tikzpicture}[trim axis right]
            \centering
            \begin{semilogxaxis}[
              scale=0.8,
              table/col sep = comma,
              xlabel = {Number of d.o.f.s at iteration~$\ell$},
              title = {$\text{CL}_\ell$}]
              \addplot+ table [x=ndofs, y=cl]
              {new_tables/efficiency_lvm_adapt_pbg.csv};
            \end{semilogxaxis}
          \end{tikzpicture}
          \caption{History of  computed loss values.}
          \label{fig:pbg-cl}
        \end{subfigure}
	\caption{Results from the adaptive algorithm applied to the PBG fiber.}
	\label{fig:pbg-all}
\end{figure}

The scale of the PBG fiber is \(L = 2.88753 \times 10^{-6}\)~m. 
The wavelength we consider is \(k = \frac{2 \pi}{8.25} \times 10^{7}\text{m}^{-1}\). 
The refractive index of the cladding and the outer region is \(n_{\rm cladding} = 1.45\), and the refractive index of the tubes is \(n_{\rm tube} = 1.8\). 
We have used a  lattice with six dielectric rods, each with a (dimensionless) radius of \(r_{\rm tube} = 0.57142857\), encircling a fiber core of radius \(r_{\rm core} = 1.0\). 
The inner radius of the PML region is \(r_{\rm outer} = 3.80952380\), and the outer radius of the PML region is \(r_{\rm PML} = 5.80952380\).

The results are  in Figure~\ref{fig:pbg-all}. The initial mesh (Figure~\ref{fig:pbg_mesh0}) is relatively coarse, although we employed smaller elements in the core (with \(h \approx 0.25\)). 
For this fiber, as seen from Figure~\ref{fig:pbg_mesh13}, the adaptive process targets refinement within the dielectric rods that surround the solid dielectric core. We see finer scale ripples within these rods being resolved in Figure~\ref{fig:pbg_phi}. 
Again, these are better seen in the longitudinal component, and are
hardly visible in the total intensity plot (in Figure~\ref{fig:pbg_intensity}). 
The convergence of the eigenvalues during the adaptive process depicted in Figure~\ref{fig:pbg-accuracy} shows a monotonic descent, unlike what was found for the anti-resonant fibers (and like the prior cases, 
the discrete cluster $\varLambda_h^{(\ell)}$ here generally had two elements).
The convergence of CL values is seen in Figure~\ref{fig:pbg-cl}.

\bigskip

To summarize the observations in this section, we have seen three examples of modern microstructured optical fibers, where our adaptive algorithm captured certain fine-scale modal features within the microstructures. 
We have also seen that without sufficient refinement in certain (perhaps unexpected) areas of the fiber geometry, the computed mode losses could be highly inaccurate.

\bigskip



\section*{Disclaimers} 
This article has been approved for public release; distribution unlimited. Public Affairs release approval {\#}AFRL-2024-1561. 
The views expressed in this article are those of the authors and do not necessarily reflect the official policy or position of the Department of the Air Force, the Department of Defense, or the U.S. government. 

\bibliographystyle{siamplain}
\bibliography{references}
\end{document}